\documentclass{article}[9pt]

\usepackage{lipsum}
\usepackage{cite}
\usepackage{url}
\usepackage{amsfonts}
\usepackage{graphicx}
\usepackage{epstopdf}
\usepackage{algorithmic}
\usepackage{amsopn}
\usepackage{dsfont}
\usepackage{enumitem}
\usepackage{epsfig}
\usepackage{tikz}
\usepackage{pgfplots}
\pgfplotsset{compat=newest}
\usetikzlibrary{plotmarks}
\usetikzlibrary{arrows.meta}
\usepackage{grffile}
\usepackage{amsmath}
\usepackage{mathtools}
\usepackage{amssymb}
\usepackage{amsthm}
\usepackage{subfig}
\usepackage{esint}
\usepackage{enumitem}

\RequirePackage[capitalize,nameinlink]{cleveref}[0.19]
\crefname{section}{section}{sections}
\crefname{subsection}{subsection}{subsections}
\Crefname{figure}{Figure}{Figures}
\crefformat{equation}{\textup{#2(#1)#3}}
\crefrangeformat{equation}{\textup{#3(#1)#4--#5(#2)#6}}
\crefmultiformat{equation}{\textup{#2(#1)#3}}{ and \textup{#2(#1)#3}}
{, \textup{#2(#1)#3}}{, and \textup{#2(#1)#3}}
\crefrangemultiformat{equation}{\textup{#3(#1)#4--#5(#2)#6}}%
{ and \textup{#3(#1)#4--#5(#2)#6}}{, \textup{#3(#1)#4--#5(#2)#6}}{, and \textup{#3(#1)#4--#5(#2)#6}}
\Crefformat{equation}{#2Equation~\textup{(#1)}#3}
\Crefrangeformat{equation}{Equations~\textup{#3(#1)#4--#5(#2)#6}}
\Crefmultiformat{equation}{Equations~\textup{#2(#1)#3}}{ and \textup{#2(#1)#3}}
{, \textup{#2(#1)#3}}{, and \textup{#2(#1)#3}}
\Crefrangemultiformat{equation}{Equations~\textup{#3(#1)#4--#5(#2)#6}}%
{ and \textup{#3(#1)#4--#5(#2)#6}}{, \textup{#3(#1)#4--#5(#2)#6}}{, and \textup{#3(#1)#4--#5(#2)#6}}
\crefdefaultlabelformat{#2\textup{#1}#3}

\allowdisplaybreaks

\DeclareMathOperator{\R}{\mathbb{R}}
\newcommand{\e}{\varepsilon}

\newcommand{\bx}{{\bf x}}

\newcommand{\by}{{\bf y}}

\newtheorem{theorem}{Theorem}
\newtheorem{corollary}[theorem]{Corollary}
\newtheorem{lemma}[theorem]{Lemma}
\newtheorem{definition}[theorem]{Definition}

\newtheorem{Remark}{Remark}

\numberwithin{theorem}{section}

\graphicspath{{./img/}}	
\newcommand{\corr}[1]{#1}
\newcommand{\ncorr}[1]{#1}

\title{An elliptic local problem with exponential decay of the resonance error for numerical homogenization}
\author{Assyr Abdulle\footnotemark[1] \and Doghonay Arjmand\footnotemark[1] \and Edoardo Paganoni\footnote{Institute of Mathematics, \'Ecole Polytechnique F\'ed\'erale de Lausanne}}
\date{}

\begin{document}
	
	\maketitle
	
	\begin{abstract}
	Numerical multiscale methods usually rely on some coupling between a macroscopic and a microscopic model. The macroscopic model is incomplete as effective quantities, such as the homogenized material coefficients or fluxes, are missing in the model. These effective data need to be computed by running local microscale simulations followed by a local averaging of the microscopic information. Motivated by the classical homogenization theory, it is a common practice to use local elliptic cell problems for computing the missing homogenized coefficients in the macro model. Such a consideration results in a first order error $O(\e/\delta)$, where $\e$ represents the wavelength of the microscale variations and $\delta$ is the size of the microscopic simulation boxes. This error\corr{, called ``resonance error'',} originates from the boundary conditions used in the micro-problem and typically dominates all other errors in a multiscale numerical method. \corr{Optimal decay of the resonance error remains an open problem,} although several interesting approaches reducing the effect of the boundary have been proposed over the last two decades. In this paper, as an attempt to resolve this \corr{problem, we propose} a computationally efficient, fully elliptic approach with \corr{exponential decay of the resonance error}.  
\end{abstract}

\textbf{Key words.}
multiscale methods, homogenization, resonance error

\textbf{AMS subject classification.}
35B27, 65L12, 74Q10

\section{Introduction}
This paper concerns the numerical homogenization of multiscale elliptic partial differential equations (PDEs) of the form
\begin{equation}\label{Eqn_Main_Multiscale_PDE}
\begin{aligned} 
-\nabla \cdot \left( a^{\e}(\bx) \nabla u^{\e}(\bx) \right) &= f(\bx) & \quad & \text{in } \Omega\,\subset \R^d \\
u^{\e}(\bx) &= 0 & \quad & \text{on } \partial\Omega,
\end{aligned} 
\end{equation}
where $a^{\e}$ characterises a microscopically non-homogeneous medium, which has small scale variations of size $\e \ll |\Omega| = O(1)$. The multiscale elliptic PDE \cref{Eqn_Main_Multiscale_PDE} is chosen to simplify the exposition, but the discussion is equally valid also for homogenization problems of parabolic and second-order hyperbolic types. Approximating the solution $u^{\e}$ via a direct numerical simulation is prohibitively expensive as accurate approximations require resolutions down to the finest scales of the problem. Sub-linear scaling multiscale numerical methods can be designed at the expense of targeting only a local average behaviour of the full solution $u^{\e}$. This local average behaviour is associated with the homogenized limit as $\e \to 0$, where the following homogenized PDE describes the coarse scale response of the system.
\begin{equation}\label{Eqn_Main_Homogenized_PDE}
\begin{aligned} 
-\nabla \cdot \left( a^{0}(\bx) \nabla u^{0}(\bx) \right) &= f(\bx) & \quad & \text{in } \Omega \\
u^{0}(\bx) &= 0 & \quad & \text{on } \partial\Omega.
\end{aligned}
\end{equation}
The existence of the homogenized solution $u^{0}$ is guaranteed by the homogenization theory, see e.g., \cite{Bensoussan_Lions_Papanicolaou,Jikov_Kozlov_Oleinik,Pavliotis_Stuart} for a few well-known monographs, in which the theoretical details of the subject are discussed. Here the homogenized coefficient $a^{0}$ is no longer oscillatory, and the problem \cref{Eqn_Main_Homogenized_PDE} may be approximated by a standard numerical method once $a^{0}$ is determined. Note that, apart from computing the homogenized solution $u^{0}$, the very goal of determining the homogenized coefficient $a^{0}$ is also important and practically relevant in many applied disciplines, e.g.,\ mechanics and material sciences.   

Explicit representations for the homogenized coefficient $a^{0}$ exist only in limited cases of interests, namely for periodic or stationary ergodic random coefficients. For example, if the coefficient $a^{\e}(\bx) = a(\bx/\e)$, where $a$ is a $K:=(-\frac{1}{2},\frac{1}{2})^{d}$-periodic function, then the homogenized coefficient is given by\footnote{Einstein summation is used in this formula. Namely, the repeated index $k$ is to mean summation over $k=1,\ldots,d$.} 
\begin{equation} \label{Eqn_Exact_Hom_Coeff}
e_i\cdot  a^{0} e_j = \fint_{K} \Big( a_{ij}(\bx ) + a_{ik}(\bx) \partial_{x_k} \chi^{j}(\bx) \Big) \; d\bx, 
\end{equation}
where $\{ e_i \}_{i=1}^{d}$ are the canonical basis vectors in $\mathbb{R}^d$, and $\{ \chi^{j} \}_{j=1}^{d}$ are the solutions of the following cell-problems posed over the unit cube $K$:
\begin{equation}\label{Eqn_Periodic_Elliptic_CellProblem}
\begin{aligned} 
-\nabla \cdot \left( a(\bx) \nabla \chi^{j}(\bx) \right)  &= \nabla \cdot a e_j, \quad \text{ in } K, \\
\chi^{j}(\bx) \quad & \text{ is } K\text{-periodic}.  
\end{aligned}
\end{equation}
From a practical point of view, the heterogeneous coefficient $a^{\e}$ is often not fully periodic and it includes more complex non-periodic variations, for which the formula \eqref{Eqn_Exact_Hom_Coeff} would either break down or be simply inaccurate. This has triggered the birth and development of a number of multiscale methodologies which target the coarse scale behaviour of the solution $u^{\e}$, without assuming an a priori knowledge about the homogenized coefficient $a^{0}$ or the precise nature of $a^{\e}$. 

Two frameworks, which address the numerical homogenization problems among other multiscale and multiphysics problems, are the Heterogeneous Multiscale Method (HMM) \cite{Abdulle_E_Engquist_1,E_Engquist_1}, and the equation free approaches \cite{Kevrikidis_etal}, which rely on a micro-macro coupling to approximate the homogenized solution $u^{0}$. In these two approaches, the scale separation in \eqref{Eqn_Main_Multiscale_PDE}, $\e \ll |\Omega|$, is exploited to design sub-linear scaling algorithms, in which unknown homogenized quantities are computed by upscaling microscopic information obtained by running local simulations over microscopic boxes of size $\delta^d$, where $\delta  = O(\e)$. A pre-selected macro model is supplied with these microscopic information and an approximation to $u^{0}$ is computed by  inverting/evolving the macroscopic problem. Other approaches, including the Multiscale Finite Element Method (MsFEM) \cite{Hou_Wu} and the Localized Orthogonal Decomposition (LOD) \cite{Henning_Malqvist_2014, Malqvist_Peterseim_2014}, address the problem in a slightly different way and aim at directly computing the oscillatory response $u^{\e}$. 
A common ingredient of these methods is the need for solving a set of local microscale problems which would then imply artificial boundary conditions on the boundaries of the microscopic domains. These boundary conditions degrade the overall accuracy of multiscale methods and improved methodologies with reduced boundary errors are needed. \ncorr{The issue of localization is also relevant for other classes of multiscale methods suited for problems  with continuum of scales and rough coefficients with high contrast. Advances in this direction include Gamblets \cite{Owhadi_2016,Owhadi_Zhang_3}, flux-transfer transformation \cite{Owhadi_Zhang_2}, polyharmonic homogenization \cite{Owhadi_Zhang_Berlyand}, and variational multiscale (VMS) based methods such as \cite{Peterseim_Verfurth,Kornhuber_Peterseim_Yserentant}. As these methods do not exploit the scale separation (for problems with scale separation), they become computationally more expensive (in comparison to HMM or equation free approaches) when applied to problems with scale separation. In what follows, we return back to our discussion on problems with scale separation and focus on coefficients with moderate contrast ratios, which are very typical in many applications, e.g., structural mechanics.} To put the discussion in a mathematical framework, we give a motivation of the boundary error by considering an example of purely periodic tensors. However, it should be kept in mind that this error is present also for more complicated non-periodic coefficients. 
\subsection{Motivation - source of the boundary error}
The generality of the multiscale algorithms, such as HMM, originates from the fact they do not suffer from the structural assumptions (other than the scale separation), that the classical analytical homogenization theory uses to derive formulas for the homogenized tensor, e.g.\ note the periodicity requirement in the formula \eqref{Eqn_Periodic_Elliptic_CellProblem}. When the period of the coefficient is not known or the medium is non-periodic, e.g., a random stationary ergodic medium, a widely used approach is to pose the cell problem \eqref{Eqn_Periodic_Elliptic_CellProblem} over a larger computational domain, say $K_{R}$ with $R>1$, and compute the approximate homogenized coefficient by an averaging over $K_R$, i.e.,  
\begin{equation} \label{Eqn_Elliptic_Approx_Hom_Coeff}
e_i \cdot a^{0}_{R} e_j =  \fint_{K_{R}} \Big( a_{ij}(\bx ) + a_{ik}(\bx) \partial_{x_k} \chi_{R}^{j}(\bx) \Big) \; d\bx, 
\end{equation}
where
\begin{equation}\label{Eqn_Approx_Periodic_Elliptic_CellProblem}
\begin{aligned} 
-\nabla \cdot \left( a(\bx)  \nabla \chi_{R}^{j}(\bx) \right)  &= \nabla \cdot a e_j, \quad \text{ in } K_{R}:=(-R/2,R/2)^d, \\
\chi_R^{j}(\bx) \quad & \text{ is } \text{ periodic in } K_{R}. 
\end{aligned}
\end{equation}
Assuming for a moment that $a$ is $K$-periodic, and $R$ is an integer, it readily follows by periodically extending $\chi$ to $K_R$ that $\chi_R = \chi$, and hence $a^{0}_R = a^{0}$. In general, when $R$ is non-integer, there is a mismatch between the values of $\chi_R$ and $\chi$ on the boundary  $\partial K_{R}$ (the so-called resonance or cell-boundary error), which yields \cite{E_Ming_Zhang}
\begin{equation}\label{eq_conv_standard_modelling_error}
\| a^{0}_{R}  -  a^{0} \|_{F} \leq C \dfrac{1}{R},
\end{equation}
where $\|\cdot\|_F$ is the Frobenius norm.
This first order resonance error, first mentioned in the context of multiscale finite element methods in \cite{Hou_Wu}, is a common drawback of modern multiscale methods as very large values for $R$ are needed to bring down this error to practical degrees of interests. For moderate values of $R$, say $R \approx 10$, this error will dominate all other errors in typical multiscale algorithms and deteriorate the overall accuracy. Hence, more efficient strategies, leading to high order rates in $1/R$, are needed. For similar results in random media see also \cite{Bourget_Piatnitski_04}. 
\subsection{Existing approaches to reduce the resonance error}
In the past, several approaches have been proposed to reduce the boundary error. The effect of using different boundary conditions (BCs), e.g.,\ Dirichlet, Neumann, or periodic BCs for the cell-problem \eqref{Eqn_Approx_Periodic_Elliptic_CellProblem} is studied in \cite{Yue_E_2007}. It is found that using different BCs does not improve the first order convergence rate in general, but periodic BCs result in a smaller prefactor compared to the Dirichlet and Neumann counterparts. There are other promising approaches which are based on modifying the cell-problem \eqref{Eqn_Approx_Periodic_Elliptic_CellProblem} so that the effect of the boundary is reduced over the interior of the domain $K_{R}$, while still retaining a good approximation of the homogenized coefficient $a^{0}$. In \cite{Blanc_LeBris_2010} a filtered cell-problem together with an integral constraint for the gradient of the cell solution is used to obtain second order convergence rates in $\frac{1}{R}$. In \cite{Gloria_2011}, a zero-th order term is added to the cell-problem so that the Green's function of the modified problem decays exponentially and the effect of mismatching BCs is significantly reduced in the interior of the domain $K_{R}$. The asymptotic convergence rate for this strategy is fourth order, but large values for $R$, e.g.,\ $R\approx 100$ are needed to observe this rate in simulations \cite{Gloria_2011}. The (fixed) convergence rate of this approach can be improved even to higher (arbitrary) orders by a Richardson extrapolation, \cite{Gloria_2012,Gloria_Habibi_2016} \corr{at the cost of iteratively solving the microscale problem}. 
\corr{Using Richardson extrapolation for improving the convergence rate of the method analysed in this paper would not be useful, as we can reach arbitrary rates of convergence without any additional cost, except the one of solving the corrector equations.}
Another approach leading to arbitrarily high orders in $\frac{1}{R}$ is proposed in \cite{Arjmand_Runborg_2016,Arjmand_Stohrer_2016}, where a second-order wave equation is used instead of the cell-problem \eqref{Eqn_Approx_Periodic_Elliptic_CellProblem}. Due to the finite speed of propagation of waves, the errors committed on the boundary of the proposed cell problems do not influence the interior solution if the computational domain is chosen sufficiently large. Although this wave approach results in a removal of the boundary error, there are a few computational challenges with this method: $i)$ the spatial domain size increases linearly with the wave speed, $ii)$ the solution of the wave equation depends on time, and therefore additional degrees of freedom are needed to approximate the cell-solution, $iii)$ practically, accurate approximations of solutions of the wave equation require high resolutions per-wavelength, which makes the method less efficient (when compared to solving an elliptic cell-problem). 

The shortcomings of the existing approaches motivate the need for designing better/alternative methodologies with improved convergence rates. In the current study, we propose an elliptic cell-problem which has exponentially decaying boundary errors. The strategy is linked to parabolic cell-problems, which was proposed in \cite{Abdulle_Arjmand_Paganoni_18a} and analysed in \cite{Abdulle_Arjmand_Paganoni_18b}, \ncorr{see also \cite{Mourrat_2019,Gloria_Neukamm_Otto_15}}. It is shown, in \cite{Abdulle_Arjmand_Paganoni_18b}, that using parabolic cell-problems one can approximate the homogenized coefficient $a^{0}$ with exponential accuracies. From a computational point of view, \ncorr{the parabolic approach requires efficient stiff time-stepping methods, which can be a numerical challenge}, \ncorr{see the numerical results section for a cost vs tolerance comparison with the parabolic approach}. The strategy proposed in this paper aims at bypassing this problem by exploiting the properties of the parabolic cell problems in \cite{Abdulle_Arjmand_Paganoni_18b} and reformulating elliptic cell-problems that mimic the behaviour of the parabolic cell-solutions with similar (but not the same) exponentially decaying convergence rates. Krylov subspace iteration is exploited to make the computational cost comparable to the cost of solving the elliptic PDE \eqref{Eqn_Approx_Periodic_Elliptic_CellProblem}. The analysis in this paper, whose results were announced in \cite{Abdulle_Arjmand_Paganoni_18a}, is done for periodic coefficients but the method itself is not limited by such a structural assumption. \ncorr{The convergence analysis for non-periodic coefficients is not addressed in the present work, however numerical results for quasi-periodic and random media are depicted to investigate the decay rates as $R$ gets larger (see Examples 3 and 4 in \Cref{Sec_Numerical_Results}).}
\subsection{Notations and definitions}
Throughout the exposition, we will use the following notations:
\begin{itemize}
\item The Sobolev space $W^{k,p}(\Omega)$ is defined as 
$$W^{k,p}(\Omega):= \{ f: D^{\gamma} f \in L^{p}(\Omega) \text{ for all multi-index } \gamma \text{ with } |\gamma| \leq k\}.$$
The norm of a function $f \in W^{k,p}(\Omega)$ is given by 
\begin{align*}
\| f \|_{W^{k,p}(\Omega)} := \begin{cases} \left( \sum_{|\gamma| \leq k} \int_{\Omega} |D^{\gamma} f(\bx)|^{p} \; d\bx \right)^{1/p} & (1 \leq p < \infty) \\
\sum_{|\gamma| \leq k} \text{ess sup}_{\Omega} |D^{\gamma} f| & (p=\infty).
\end{cases}
\end{align*}
\item The space $H^1_0(\Omega)$ is the closure in the $W^{1,2}$-norm of $C_c^{\infty}(\Omega)$, the space of infinitely differentiable functions with compact support in $\Omega$. The norm associated with $H^1_{0}(\Omega)$ is 
	\begin{equation*}
	\| f \|^2_{H_{0}^1(\Omega)} := \| f \|^2_{L^2(\Omega)} + \| \nabla  f \|^2_{L^2(\Omega)} ,
	\end{equation*}
\corr{but an equivalent norm is
	\begin{equation*}
	\| f \|_{H_{0}^1(\Omega)} := \| \nabla  f \|_{L^2(\Omega)}.
	\end{equation*}
We will use this second notation for the $H_{0}^1$-norm.
	}
\item We use the notation $\langle f, g \rangle_{L^2(\Omega)} := \int_{\Omega} f g \; d\bx$ to denote the $L^2$ inner product over $\Omega$.  
\item The space $H_{\text{div}}$ is 
$$H_{\text{div}}(\Omega) := \{f: f \in [L^2(\Omega)]^{d} \text{ and } \nabla \cdot f \in L^{2}(\Omega)\}.$$
The norm associated with $H_{\text{div}}$ is 
\begin{align*}
\| f \|_{H_{div}(\Omega)}^{2} := \| f \|_{L^2(\Omega)}^{2} +  \| \nabla \cdot f \|_{L^2(\Omega)}^{2}. 
\end{align*}
\item Cubes in $\mathbb{R}^{d}$ are denoted by $K_{L}:=(-L/2,L/2)^{d}$. In particular, $K$ is the unit cube. 
\item  \corr{The space $W^1_{per}(K)$ is defined as the closure of 
$$\left\{ f \in C^{\infty}_{per}(K): \int_{K} f = 0 \right\}$$
with respect to the $W^{1,2}$-norm. Thanks to the Poincar\'e-Wirtinger inequality, we can also have the following equivalent norm
\begin{equation*}
\| f \|_{W^1_{per}(K)} :=  \| \nabla  f \|_{L^2(K)}. 
\end{equation*}}
\item Let $f$ belong to the Bochner space $ L^{p}(0,T; X)$, where $X$ is a Banach space. Then the norm associated with this space is defined as 
$$
\| f \|_{L^{p}(0,T; X)}:= \left( \int_{0}^{T} \|  f \|_{X}^{p} \; dt \right)^{\frac{1}{p}}. 
$$
\item By writing $C$, we mean a generic constant independent of $R,L,T,N$ which may change in every subsequent occurrence.  
\item Boldface letters \corr{in arguments of functions} are to distinguish functions in multi-dimensions, e.g., $f(\bx)$ is to mean a function of several variable ($\bx \in \mathbb{R}^d, d \geq 2$), while $f(x)$ will be a function of one variable ($ x \in \mathbb{R}$).  
\item We will use the notation $\fint_{D} f(\bx) \; d\bx$ to denote the average $\frac{1}{D}\int_{D}f(\bx) \; d\bx$ over a domain D.
\end{itemize}
\begin{definition}[Filter as in \cite{Gloria_2011}]\label{Def_Filter} We say that a function $$\mu: [-1/2,1/2] \to \mathbb{R}^{+}$$ belongs to the space $\mathbb{K}^{q}$ with $q > 0$ if
\begin{itemize}
\item[i)] $\mu \in C^{q}([-1/2,1/2]) \cap W^{q+1,\infty}((-1/2,1/2))$
\item[ii)] $\int_{-1/2}^{1/2}  \mu(x) \; dx  =1, $
\item[iii)] $\mu^{k}(-1/2) = \mu^{k}(1/2) =  0$ for all $k \in \{ 0,\ldots,q-1 \}$.
\end{itemize}
In multi-dimensions a $q$-th order filter $\mu_{L}: K_{L} \to \mathbb{R}^{+}$ with $L>0$ is defined by  
\begin{align*}
\mu_L(\bx) := L^{-d} \prod_{i=1}^{d} \mu\left(\frac{x_i}{L}\right),
\end{align*}
where $\mu$ is a one dimensional $q$-th order filter and $\bx = (x_1,x_2,\ldots,x_d) \in \mathbb{R}^d$. \corr{In this case, we will say that $\mu_L \in \mathbb{K}^{q}(K_L)$}. Note that filters $\mu_L$ are considered extended to 0 outside of $K_L$.
\end{definition}
Filters have the property of approximating the average of periodic functions with arbitrary rate of accuracy, as state in the following \cref{Lem_Filtering} (see \cite{Gloria_2011} for a proof).
\begin{lemma}\label{Lem_Filtering}
		Let $\mu_L \in \mathbb{K}^{q}(K_L)$. Then, for any $K$-periodic function $f\in L^p(K)$ with $1 < p\le 2$, we have 
		\begin{equation*} 
		\left| \int_{K_L}  f(\bx) \mu_L(\bx) \, d\bx - \int_K f(\bx) \,d\bx \right| \le C_{\mu} \left\|f\right\|_{L^{p}(K)} L^{-(q+1)},
		\end{equation*}
		\ncorr{where $C_{\mu}$ is the Lipschitz constant of $q$-th derivative of the filter $\mu$.}
\end{lemma}
\begin{definition}\label{Def_MSpace} We say that $a \in \mathcal{M}(\alpha,\beta,\Omega)$ if $a_{ij} = a_{ji}$, $a \in [L^{\infty}(\Omega)]^{d\times d}$ and there are constants $0<\alpha \leq \beta$ such that
\begin{align*}
\alpha |\zeta|^2 \leq a(\bx) \zeta \cdot \zeta \leq \beta |\zeta|^2, \quad \text{ for a.e. } \bx \in \Omega, \forall   \zeta \in \mathbb{R}^d.
\end{align*} 
We write $a \in \mathcal{M}_{\text{per}}(\alpha,\beta,\Omega)$ if in addition $a$ is a $\Omega$-periodic function. 
\end{definition}
This paper is structured as follows. In \cref{Sec_Fully_Elliptic_Approach}, we propose a fully elliptic approach with exponentially decaying resonance errors. In \cref{Sec_Main_Results}, we present the main results and provide an analysis of the exponential convergence rates. 
A numerical strategy, based on the Arnoldi decomposition, to approximate the solution of the proposed cell problem is discussed in \cref{Sec_Arnoldi}.
Finally, \cref{Sec_Numerical_Results} includes numerical tests supporting the theoretical findings.    
\section{A modified elliptic approach}\label{Sec_Fully_Elliptic_Approach}
The approach proposed here is based on adding a correction term to the elliptic cell problem \eqref{Eqn_Approx_Periodic_Elliptic_CellProblem} so that the effect of the boundary values are significantly reduced. For $j=1,\ldots,d$, the new cell-problems read as
\begin{equation}\label{Eqn_New_CellProblem}
\begin{aligned} 
-\nabla \cdot \left( a(\bx) \nabla \chi_{T,R}^{j}(\bx) \right)  &= g^{j}(\bx) - [e^{-A T} g^{j}](\bx) \quad \text{ in } K_{R} \\
\chi_{T,R}^{j}(\bx) &= 0 \quad \text{ on } \partial K_{R}.  
\end{aligned}
\end{equation}
Here $A$ and $g^{j}$ are defined as 
\begin{align*}
A:= -\nabla \cdot \left( a \nabla  \right), \text{ and } g^{j} := \nabla \cdot (a e_j),
\end{align*}
where we implicitly assume that $A$ is equipped with the homogeneous Dirichlet boundary conditions. Note that $a^{\e}(\bx):=a(\bx/{\e})$ does not need to be periodic, but the periodicity of $a$ will be assumed later for the analysis. The evolution operator $e^{-A T}$ is the semigroup generated by the operator $-A$, i.e., $e^{-AT} g^j$ is the solution at time $T$ of the corresponding parabolic problem ($\partial_t u  - \nabla \cdot \left( a(\bx) \nabla u \right)= 0$) over $K_R$ with homogeneous Dirichlet boundary conditions and initial data $g^j$. In a nutshell, we will present algorithms based on spectral truncation as well as a Krylov subspace iteration to approximate the correction term  $e^{-AT} g^j$ without solving the parabolic PDE. The homogenized coefficient can then be approximated by 
\begin{align} \label{Eqn_New_Approx_Hom_Coeff}
e_i\cdot a^{0}_{T,R,L} e_j = \int_{K_{L}}  \Big( a_{ij}(\bx ) + a_{ik}(\bx) \partial_{x_k} \chi_{T,R}^{j}(\bx) \Big) \mu_L(\bx) \; d\bx.
\end{align}
To get optimal rates, the parameters $T$ \ncorr{and $L$} should be chosen as a function of $R$ and the coercivity and boundedness constants $\alpha, \beta$. The precise choice will be clarified later in the sequel. The choice of the Dirichlet boundary conditions is only for theoretical purposes, but a periodic BC similar to \eqref{Eqn_Elliptic_Approx_Hom_Coeff} would also work equally well in practice. The main difference between \eqref{Eqn_New_CellProblem} and the cell-problem \eqref{Eqn_Elliptic_Approx_Hom_Coeff} is the addition of a correction term of the form $e^{-A T} g^{j}$, which is crucial to obtain exponentially decaying convergence rates for the boundary error. One other important component of the proposed method is the presence of the filter $\mu_L$ in \eqref{Eqn_New_Approx_Hom_Coeff}. Clearly, one can see from formula \eqref{Eqn_New_Approx_Hom_Coeff} that the computation of homogenized coefficients is associated with the averages of oscillatory functions. Such filters are typically used to accelerate the convergence even for general non-periodic variations, and therefore their presence is vital for improving the accuracy in the present context too.   

\ncorr{
\textbf{Integrating the cell-problem \eqref{Eqn_New_CellProblem} in FE-HMM.} Incorporating the cell-problem \eqref{Eqn_New_CellProblem} into existing finite element/finite difference heterogeneous multiscale methods (FE/FD-HMM) follows the by now classical methodology, see \cite{Abdulle_05,Abdulle_09}. For example, the main component of FE-HMM is the approximation of the bilinear form $B^{0}(u^{0},v) = \int_{\Omega} a^{0}(\bx) \nabla u^{0} \cdot \nabla v(\bx)\; d\bx$, where $a^{0}(\bx)$ is an unknown homogenized coefficient, and $u^{0}$ is the homogenized solution. If $a^{0}(\bx)$ were to be known, then the bilinear form $B^{0}(u^{0},v)$ could be approximated by summing the following contributions over all triangles $K$ in a macro triangulation $\mathcal{T}_H$ as $B^{0}_H(u_H^{0},v) = \sum_{K \in \mathcal{T}_{H} } \omega_{K,J} a^{0}(\bx_{K,J}) \nabla u^{0}_{H}(\bx_{K,J}) \cdot \nabla v_{H}(\bx_{K,J})$, where $\omega_{K,J}$ are appropriate weight functions and $\bx_{K,J}$ are the corresponding quadrature nodes in the macroscopic mesh. In order to compute $a^{0}(\bx_{K,J})$, one can then rescale the problem \eqref{Eqn_New_CellProblem}, and the upscaling formula \eqref{Eqn_New_Approx_Hom_Coeff} to $ K_{\varepsilon R}$, and center it at $\bx_{K,J}$. The rescaled problem reads as  
\begin{equation*}
\begin{aligned} 
-\nabla \cdot \left( a^{\varepsilon}(\bx) \nabla \chi_{\varepsilon^2 T,\varepsilon R}^{j}(\bx) \right)  &= g^{j,\varepsilon}(\bx) - [e^{-A^{\varepsilon} \varepsilon^2 T} g^{j,\varepsilon}](\bx) \quad \text{ in } \bx_{K,J} +  K_{\varepsilon R} \\
\chi_{\varepsilon^2 T,\varepsilon R}^{j}(\bx) &= 0 \quad \text{ on } \partial K_{\varepsilon R},  
\end{aligned}
\end{equation*}
where $A^{\varepsilon} :=-\nabla \cdot \left( a^{\varepsilon}(\bx) \nabla\right), $ and $g^{j,\varepsilon} = -\nabla \cdot a^{\varepsilon}(\bx) e_j$.
Moreover, the upscaling formula will be given by 
\begin{align*}
e_i\cdot a^{0}(\bx_{K,J}) e_j = \int_{\bx_{K,J} + K_{ \varepsilon L}}  \Big( a_{ij}^{\varepsilon}(\bx ) + a_{ik}^{\varepsilon}(\bx) \partial_{x_k} \chi_{\varepsilon^2 T,\varepsilon R}^{j}(\bx) \Big) \mu_{\varepsilon L}(\bx) \; d\bx.
\end{align*}
By computing the above cell problem using a micro triangulation $\tau_h$, an approximation $a^{0}_h(\bx_{K,J})$ of $a^{0}(\bx_{K,J})$ is obtained and the overall algorithm relying on a macro and micro triangulation can be readily implemented. The overall computational cost of such multiscale algorithms with respect to the macro-micro mesh sizes are well-known, see e.g., \cite{Abdulle_05,Abdulle_09}, but in this paper, we also discuss the computational complexity of the cell-problem \eqref{Eqn_New_CellProblem} in Section \ref{SubSec_CompCost}.
}

\begin{lemma}[Well-posedness of \eqref{Eqn_New_CellProblem}] Let $a \in \mathcal{M}(\alpha,\beta,K_{R})$ and $ae_j \in H_{\text{div}}(K_{R})$. Then there exists a unique weak solution of \eqref{Eqn_New_CellProblem} in $H_{0}^1(K_{R})$ satisfying the estimate
\begin{equation}\label{Estimate_ChiTR}
\| \chi_{T,R}^{j} \|_{H_0^{1}(K_{R})} \leq C  \|  a e_{j}  \|_{H_{div}(K_{R})}, 
\end{equation}
where $C^2=\frac{1 + C_p^2}{\alpha^2}$, $C_p$ is the Poincar\'e constant and  $\alpha$ is the coercivity constant.
\end{lemma}

\begin{proof}
The operator $A$ generates the contraction semigroup (i.e. a family of bounded linear operators) $e^{-AT}: L^{2}(K_{R}) \to L^{2}(K_{R})$ which satisfies
$$\| e^{-A T} \|_{L^2(K_{R}) \to L^{2}(K_R)} \leq 1.$$
Since $g^{j}:=\nabla \cdot a e_j \in L^{2}(K_{R})$, we have
\begin{equation*}
\| e^{-AT} g^{j} \|_{L^2(K_{R})} \leq \| g^{j} \|_{L^2(K_R)}\quad \forall T\ge 0.
\end{equation*}
To see this, we recall that the eigenvalues of $A$ are positive and satisfy  
\begin{align*}
	0<\lambda_0 \leq \lambda_1 \leq \lambda_2, \ldots,
\end{align*}
and the eigenfunctions $\{ \varphi_k \}_{k=0}^{\infty}$ form an orthonormal basis for $L^2(K_{R})$. Next, we write
\begin{align*}
e^{-A T} g^{j}(\bx)  = \sum_{k=0}^{\infty} e^{-\lambda_k T} g^{j}_{k} \varphi_{k}(\bx),\quad \text{with } g^{j}_{k} := \langle g^{j} , \varphi_{k} \rangle_{L^2(K_{R})}.
\end{align*}
Since the eigenfunctions are orthonormal in $L^2(K_{R})$, it follows that 
\begin{equation}\label{L2_est_gT}
\| e^{-A T} g^{j} \|_{L^2(K_{R})}^2  = \sum_{k=0}^{\infty} e^{-2 \lambda_k T} |g^{j}_k|^2 \leq  \sum_{k=0}^{\infty} |g^{j}_k|^2  = \| g^{j} \|^{2}_{L^2(K_{R})}.  
\end{equation}
The Lax-Milgram theorem guarantees the existence and uniqueness of $\chi_{T,R}$ in the space $H_{0}^1(K_{R})$. By uniform ellipticity of the coefficients $a(\bx)$ and H\"older inequality we derive that
\begin{equation*}
\alpha \| \nabla \chi^j_{T,R} \|_{L^{2}(K_{R})}^2 
\le \| ae_j\|_{L^2(K_R)} \|\nabla\chi^j_{T,R}\|_{L^2(K_R)} + \|e^{-A T} g^{j}\|_{L^2(K_R)} \|\chi^j_{T,R}\|_{L^2(K_R)},
\end{equation*}
that, by application of \eqref{L2_est_gT} and the Poincar\'e inequality for $\chi_{T,R}$ and Young inequality, leads to the final bound 
\begin{equation*}
\corr{\| \chi^j_{T,R} \|_{H_0^{1}(K_{R})} \leq \sqrt{\frac{1 + C_p^2}{\alpha^2}}\|  ae_j  \|_{H^1_{div}(K_{R})}, }
\end{equation*}
where $\alpha$ is the ellipticity constant and $C_p$ is the Poincar\'e constant. 
\end{proof}

\begin{Remark}
Note that periodicity of $a$ is not necessary for the well-posedness of $\chi_{T,R}^{j}$. 
	\ncorr{The constant in the estimate \eqref{Estimate_ChiTR} depends on the Poincar\'e constant hence, implicitly, on the size of the domain $R$. The dependency on $R$ is not stated explicitly because the convergence result stated in \cref{Thm_Main_Result} does not rely on \eqref{Estimate_ChiTR}, but only on the well-posedness of \eqref{Estimate_ChiTR}.}
\end{Remark}

\subsection{Relation with parabolic cell problems} 
The use of parabolic cell problems, that results in exponential decay of the boundary error, has been recently proposed in \cite{Abdulle_Arjmand_Paganoni_18a} and analysed in \cite{Abdulle_Arjmand_Paganoni_18b}. 
\corr{This theory is based on an idea developed earlier in \cite{Mourrat_2019}.}
The main idea behind this exponential decay is that, over sufficiently small time frames, the solutions of parabolic PDEs over a bounded domain ``do not feel the boundary''. In fact, the proposed elliptic cell problems \eqref{Eqn_New_CellProblem} are closely related to the solution of parabolic PDEs \corr{and this relation will be used in the subsequent analysis}.
In this section, we present a theorem which states that the time integration of parabolic cell problems solve the elliptic cell problem \eqref{Eqn_New_CellProblem}. 

\begin{theorem}\label{Thm_Time_Integration}
Assume that $a \in \mathcal{M}(\alpha,\beta,K_{R})$ and $ae^{j} \in H_{div}(K_{R})$ and let $u^{j}$ be the solution of the following parabolic PDE
\begin{equation}\label{Eqn_Parabolic_CellProblem}
\begin{aligned} 
\partial_{t} u^{j}(t,\bx) &- \nabla \cdot \left( a(\bx) \nabla u^{j}(t,\bx) \right) = 0, \quad \text{in }  K_{R} \times (0,T],\\
u^{j}(t,\bx) &= 0, \quad \text{on }  \partial K_{R} \times (0,T]  \\
u^{j}(0,\bx) &= g^{j}(\bx), \quad \text{in }  K_{R}, 
\end{aligned}
\end{equation}
where $g^{j}(\bx):=\nabla \cdot \left( a(\bx) e_{j}\right)$. Then, the time integral
\begin{equation*}
\chi^{j}_{T,R}(\bx)  = \int_{0}^{T} u^{j}(t,\bx) \; dt
\end{equation*}
solves the PDE
\begin{equation}\label{Eqn_CellProblem}
\begin{aligned} 
-\nabla \cdot \left( a(\bx) \nabla \chi^{j}_{T,R}(\bx) \right)  &= g^{j}(\bx) - [e^{-A T} g^{j}](\bx), &\quad &\text{in } K_{R} \\
\chi^{j}_{T,R}(\bx) &= 0 & \quad & \text{on } \partial K_{R}, 
\end{aligned}
\end{equation}
where $A:= -\nabla \cdot \left( a \nabla  \right)$.
\end{theorem}
\begin{proof}
Let $\lbrace\varphi_{R,k}\rbrace_{k=0}^{\infty}$ be the eigenfunctions of the operator $A$.
\ncorr{The subscript $R$ is used to highlight the dependency of the eigenfunctions (and eigenvalues) on the size of the domain, $R$.}
The expansion of the solution $u^{j}$ in terms of $\lbrace\varphi_{R,k}\rbrace_{k=0}^{\infty}$ gives
\begin{equation*}
u^{j}(t,\bx)  = \sum_{k=0}^{\infty} u_k^{j}(t) \varphi_{R,k}(\bx), \quad u^{j}(0,\bx)  = \sum_{k=0}^{\infty} g^j_{k} \varphi_{R,k}(\bx),
\end{equation*}
where $g^{j}_{k} := \langle g^{j} , \varphi_{k} \rangle_{L^2(K_{R})}$ and $u^{j}_{k} := \langle u^{j} , \varphi_{k} \rangle_{L^2(K_{R})}$.
Plugging this expansion into the equation \eqref{Eqn_Parabolic_CellProblem}, we obtain
\begin{equation*}
\sum_{k=0}^{\infty} \left(\dfrac{d}{dt}u_k^{j}(t) \varphi_{R,k}(\bx) + u_k^{j}(t) \lambda_{R,k} \varphi_{R,k}(\bx)\right)  = 0, \quad \text{ for } j = 1,\dots,d,
\end{equation*}
where $\lbrace\lambda_{k}\rbrace_{k=0}^{\infty}$ are the eigenvalues of $A$.
Since the eigenfunctions are orthonormal in $L^2(K_{R})$, we arrive at 
\begin{equation*}
u_k^{j}(t)  = e^{-\lambda_{R,k} t} u_k^j(0)  = e^{-\lambda_{R,k} t} g_k^{j}.
\end{equation*}
Now, integrating in time, we obtain
\begin{align*}
\chi^{j}_{T,R}(\bx) &:= \int_{0}^{T} u^{j}(t,\bx) \;dt \\
 &= \sum_{k=0}^{\infty} g_k^{j} \varphi_{j}(\bx) \int_{0}^{T}  e^{-\lambda_{R,k} t} \; dt \\
 &=  \sum_{k=0}^{\infty} \dfrac{1}{\lambda_{R,k}} g_k^{j} \varphi_{R,k}(\bx) - \sum_{k=0}^{\infty} \dfrac{1}{\lambda_{R,k}} e^{-\lambda_{R,k} T} g_k^{j} \varphi_{R,k}(\bx)
\end{align*}
Moreover, evaluating $A \chi_{T,R}(\bx)$ we obtain
\begin{align*}
	A \chi^{j}_{T,R}(\bx) &=  \sum_{k=0}^{\infty} g_k^{j} \varphi_{R,k}(\bx) - \sum_{k=0}^{\infty} e^{-\lambda_{R,k} T} g_k^{j} \varphi_{R,k}^{j}(\bx) \\
	&= g^{j}(\bx) - e^{-AT} g^{j}(\bx).
\end{align*}

\end{proof} 

\ncorr{By integrating equation \eqref{Eqn_Parabolic_CellProblem} in time from $t=0$ to $t= \infty$, and using the fact that $u^{j} \to 0$ as $t\to \infty$, it is easy to see that the correction term vanishes, and the standard elliptic cell-problems with Dirichlet boundary conditions are recovered; i.e., the cell-problem \eqref{Eqn_Approx_Periodic_Elliptic_CellProblem} but with homogeneous Dirichlet conditions.} Therefore, no improvement in terms of the overall convergence rate for the resonance error will be observed if the parameter $T$ ``is chosen too large''. Similarly, by studying the limiting equation as $T \to \infty$ of parabolic PDEs with periodic solutions over the unit cell $K$, one can recover the periodic cell-problem \eqref{Eqn_Periodic_Elliptic_CellProblem}. This result is stated as a theorem here and will be used later in the analysis. 

\begin{theorem}\label{Thm_Time_Integration_Periodic}
Assume that $a \in \mathcal{M}_{\text{per}}(\alpha,\beta,K)$ and $ae^{j} \in H_{div}(K)$. Moreover, let $v^{j}$ be the solution of the following parabolic PDE 
\begin{equation}\label{Eqn_Parabolic_Periodic_CellProblem}
\begin{aligned} 
\partial_{t} v^{j}(t,\bx) &- \nabla \cdot \left( a(\bx) \nabla v^{j}(t,\bx) \right) = 0, \quad \text{in }  K \times (0,T],\\
v^{j}(t,\bx) &\quad \text{ is periodic in } K \\
v^{j}(0,\bx) &= \nabla \cdot \left( a(\bx) e_{j}\right), \quad \text{in }  K .
\end{aligned}
\end{equation}
Then the time integral
\begin{equation*}
\chi^{j}(\bx)  = \int_{0}^{\infty} v^{j}(t,\bx) \; dt
\end{equation*}
solves the PDE
\begin{align*}
-\nabla \cdot \left( a(\bx) \nabla \chi^{j}(\bx) \right)  = \nabla \cdot \left( a(\bx) e_{j}\right), \quad \text{ in } K \\
\chi^{j}(\bx) \quad \text{ is periodic in } K. \nonumber 
\end{align*}
\end{theorem}
\corr{We now give an \textit{a-priori} estimate on $\nabla v^j$, which will be used in the subsequent analysis of the resonance error. The proof is based on the spectral properties of the periodic cell-problem.}
\begin{lemma}\label{Lem_Estimate_Periodic_Parabolic_Problem}
Let us assume that the hypothesis of \cref{Thm_Time_Integration_Periodic} hold true and let $g^{j}(\bx) := \nabla \cdot \left( a(\bx) e_{j}\right)$. Then, the solution $v^j$ of \eqref{Eqn_Parabolic_Periodic_CellProblem} satisfies 
\begin{equation}\label{Eq1_Estimate_Periodic_Parabolic_Problem}
\| \nabla v^{j} \|_{L^{2}(0,T;L^{2}(K))} \leq C_1(\alpha) \| g^{j} \|_{L^{2}(K)}.  
\end{equation}
Moreover, if $g^{j} \in W^{1}_{\text{per}}(K)$, then
\begin{equation}\label{Eq2_Estimate_Periodic_Parabolic_Problem}
\| \nabla v^{j} \|_{L^{1}(0,T;L^{2}(K))} \leq C_2(\alpha,\beta) \| \nabla g^{j} \|_{L^{2}(K)},  
\end{equation}
where \ncorr{$C_1(\alpha) = \frac{1}{2\alpha}$, and $C_2(\alpha,\beta) = \frac{d \sqrt{\beta}}{\pi^2 \alpha^{3/2}}   $ }. 
\end{lemma}

\begin{proof}
Let us define the bilinear form $B:W^{1}_{per}(K)\times W^{1}_{per}(K) \mapsto \R$ as  
\begin{equation*}
B[w,\hat{w}] := \int_K \nabla \hat{w}(\bx) \cdot a(\bx) \nabla w(\bx) \,d\bx, \quad w,\hat{w} \in W^{1}_{per}(K).
\end{equation*}
With a slight abuse of notation, we denote the eigenvalues and eigenfunctions of $B[\cdot,\cdot]$ by $\{\lambda_{k}\}_{k=0}^{\infty}$ and $\{\varphi_{k}\}_{k=0}^{\infty}$, respectively. It is well known that the sequence of eigenvalues is positive and non-decreasing, i.e.,
\begin{equation*}
0<\lambda_0 \leq \lambda_1 \leq \lambda_2, \ldots.
\end{equation*}
The eigenfunctions $\{\varphi_{k}\}_{k=0}^{\infty}$ are orthonormal in the $L^2$-sense and they satisfy:
\begin{equation*}
B[\varphi_k,u] = \lambda_k \langle \varphi_k, u \rangle_{L^2(K)} ,\quad\forall u\in W^1_{per}(K).
\end{equation*}
Since the eigenvalues form a basis of $W^1_{per}(K)$, we can write the solution $v^j$ of \cref{Eqn_Parabolic_Periodic_CellProblem} as $v^{j}(t,\bx)  = \sum_{k=0}^{\infty} v_k^{j}(t) \varphi_k(\bx) $. By coercivity of the bilinear form \corr{and the exponential decay of the components of $v^j$ in the eigenfunctions' basis, $v^j_k(t) = e^{-\lambda_kt}g^j_k$ }, we obtain 
\begin{align*}
\alpha \| \nabla v^{j}(t,\cdot) \|_{L^{2}(K )}^{2} &\leq B[v^{j},v^{j}](t) \\
&= \sum_{k,\ell=0}^{\infty} e^{-(\lambda_{k} + \lambda_{\ell}) t} g^{j}_{k} g^{j}_{\ell} B[\varphi_{k},\varphi_{\ell}] \\
&\leq  \sum_{k,\ell=0}^{\infty} e^{-(\lambda_{k} + \lambda_{\ell}) t} g^{j}_{k} g^{j}_{\ell}  \lambda_k \langle \varphi_{k},\varphi_{\ell} \rangle_{L^{2}(K )} \\ 
&= \sum_{k=0}^{\infty} e^{-2\lambda_{k} t} | g^{j}_{k} |^2 \lambda_k. 
\end{align*}
From here and Parseval identity, it follows that 
\begin{align*}
\| \nabla v^{j} \|_{L^2(0,T;L^{2}(K ))}^{2} &:= \int_{0}^{T}  \| \nabla v^{j}(t,\cdot) \|_{L^{2}(K )}^{2} \; dt \\ 
&\leq  \alpha^{-1} \sum_{k=0}^{\infty} \lambda_k  \int_{0}^{T}  e^{-2\lambda_k t} \; dt   | g^{j}_{k} |^2 \\
&\leq \dfrac{\alpha^{-1}}{2} \sum_{k=0}^{\infty} | g^{j}_{k} |^2 \\
&=  \dfrac{\alpha^{-1}}{2} \| g^{j} \|_{L^2(K )}^{2}. 
\end{align*}
To prove the bound in $L^1(0,T;L^2(K ))$, we proceed as follows:
\begin{align*}
\| \nabla v^{j} \|_{L^1(0,T;L^{2}(K ))}&:= \int_{0}^{T}  \| \nabla v^{j}(t,\cdot) \|_{L^{2}(K )} \; dt \\ 
&\leq  \alpha^{-1/2} \int_{0}^{T} \sqrt{\sum_{k=0}^{\infty}     e^{-2\lambda_k t}  \lambda_k | g^{j}_{k} |^2 } \; dt\\
&\leq  \alpha^{-1/2} \int_{0}^{T} e^{-\lambda_0 t}  \; dt \sqrt{\sum_{k=0}^{\infty}  \lambda_k | g^{j}_{k} |^2 } \\
&= \dfrac{\alpha^{-1/2}}{\lambda_0} \left( 1- e^{-\lambda_0 T}  \right)  \sqrt{ B[g^{j},g^j] } \\ 
&\leq \sqrt{\frac{\beta}{\alpha}}\dfrac{1}{\lambda_0} \| \nabla g^{j} \|_{L^2(K)}.
\end{align*}
\ncorr{Here $\lambda_0 \geq \alpha C_p(K)^{-2}$, and $C_p(K)$ is the constant of the Poincar\'e-Wirtinger inequality in $W^1_{per}(K)$, which can be bounded by $C_p(K) \leq \frac{diam(K)}{\pi} = \frac{\sqrt{d}}{\pi}$, \cite{Payne_Weinerberger_1960}. Hence, $\lambda_0 \geq \frac{\alpha\pi^2}{d}$ and the final result follows.}
\end{proof}
\section{Exponential decay of the resonance error for the modified elliptic approach}\label{Sec_Main_Results}
The main result of this article is the following theorem, which gives an error bound for the difference between the exact homogenized coefficient \eqref{Eqn_Exact_Hom_Coeff} and the approximation \eqref{Eqn_New_Approx_Hom_Coeff} for a periodic material coefficient $a$. \ncorr{For proving the following \cref{Thm_Main_Result} we assume:
\begin{subequations}
	\label{eq: hp on a} 
	\begin{equation}
		\label{eq: hp on a standard}
		a(\cdot) \in \mathcal{M}_{\text{per}}(\alpha,\beta,K), 
	\end{equation}	
	\begin{equation}
		\label{eq: hp on a regularity}
		a(\cdot)\mathbf{e}_i \in W^{1,\infty}(K_R) ,\quad\text{for }i = 1,\dots,d, 
	\end{equation}	
	\begin{equation}
		\label{eq: hp on v regularity}
		v^i\in L^p\left((0,T), W^{1,p}_{per}(K)\right),\quad\text{for }i = 1,\dots,d,
	\end{equation}	
\end{subequations}
with $p>p_0=\max\left\lbrace\frac{d+2}{2},2\right\rbrace$. 
}

\begin{figure}[ht]%
	\centering
	\subfloat[$a(\bx) = diag(a_{11}(\bx),a_{22}(\bx))$]{{\includegraphics[height=4.65cm]{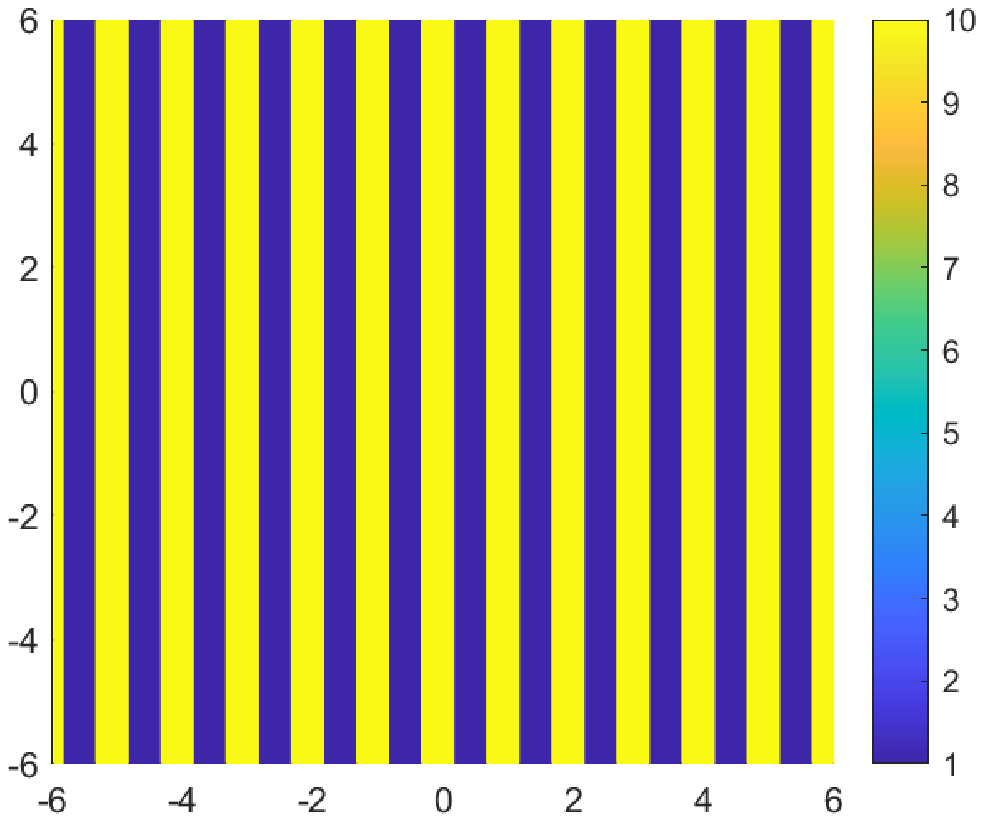} }}%
	\subfloat[The upscaling error]{{\includegraphics[height=4.65cm]{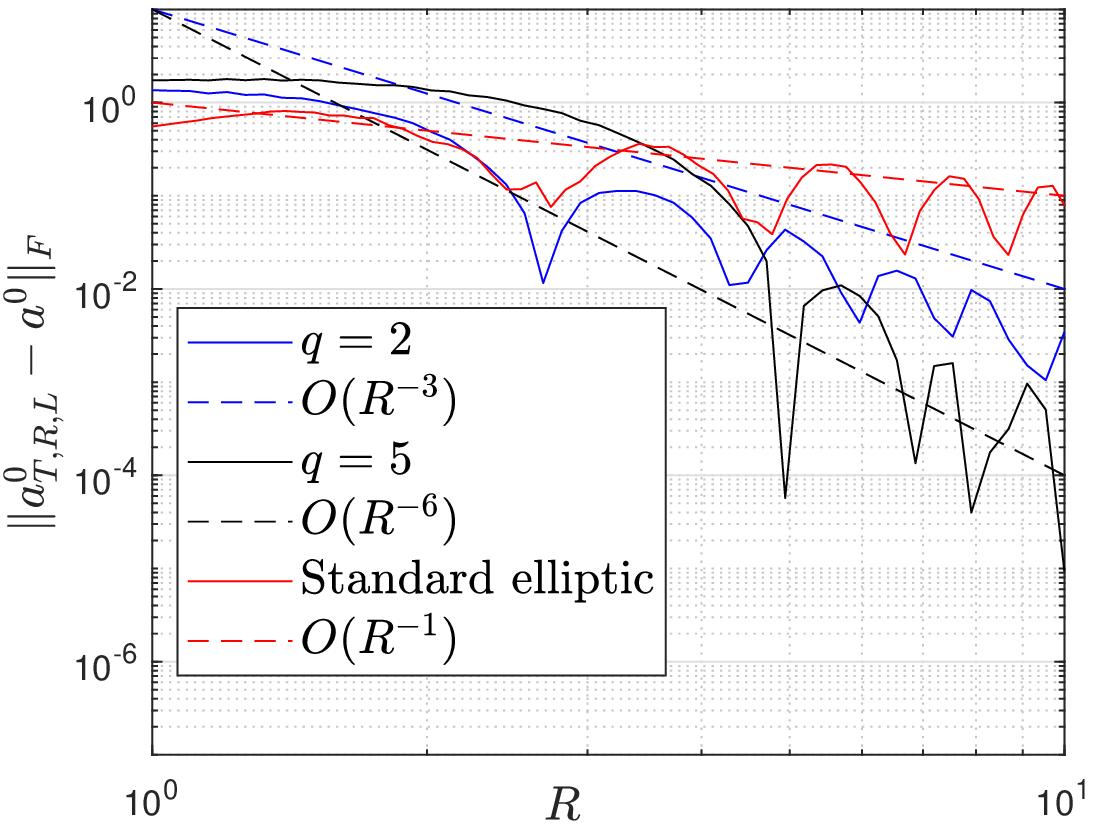}}}%
	\caption{A two dimensional periodic discontinuous medium \eqref{eqn_Coeff_Layered2D}}%
	\label{Fig_TwoD_Layered}%
\end{figure}

\ncorr{Note that the assumption \eqref{eq: hp on v regularity} has been discussed in [\cite{Abdulle_Arjmand_Paganoni_18b}, Remark 3.2]. This assumption is related to regularity of the Green's function associated with the parabolic problem \eqref{Eqn_Parabolic_Periodic_CellProblem}. We notive that for $d=1$ it is the standard estimate. For $d=2,3$, it is only slightly more stringent than the known a priori estimate $v \in  L^2\left((0,T), W^{1,2}_{per}(K)\right) $.}
\begin{theorem} \label{Thm_Main_Result}
\ncorr{
Under assumptions \eqref{eq: hp on a}, let $K_R\subset \R^d$ for $R\ge 1$,  ${\mu_L \in \mathbb{K}^{q}(K_L)}$ with $0<L<R-2$. There exists constants $c_1,c_2,C_1,C_2,C_3>0$ such that, for any $0<T<\frac{2c_2|R-L|^2}{d}$,  }
\ncorr{ 
\begin{align*}
\|  a^{0}_{T,R,L} -  a^{0} \|_{F}   &\leq \sup_{j}\left( \|\nabla \cdot a e_j \|_{L^{2}(K)} + \| \nabla \cdot a e_j \|^2_{L^{2}(K)}\right) \\ &\left( C_1 \sqrt{T} L^{-(q+1)} + C_2 e^{- c_1 T} + C_3 \left(1 + \frac{T}{|R-L|}\right) \frac{R^{d-1}}{T^{d/2-1}} e^{-c_2 \frac{\left| R-L \right|^2}{T}}\right),
\end{align*}
}
\ncorr{where $a^{0}$ and $a^{0}_{T,R,L}$ are defined in \eqref{Eqn_Exact_Hom_Coeff} and \eqref{Eqn_New_Approx_Hom_Coeff}. Moreover, $c_1 = \frac{\alpha \pi^2}{d}$, $C_1 = C_{\mu} \frac{\beta}{\alpha}$ where $C_{\mu}$ is the Lipschitz constant for the $q$-th derivative of the filter $\mu$, $C_2(d,\alpha) = \frac{d}{\alpha \pi^2}$, $C_3 = \gamma \dfrac{\beta^2}{\alpha \sqrt{\alpha c_2}}$, for a $\gamma>0$ independent of $a,R,L,T$, and $0< c_2 <\nu(\beta,d)$, where $\nu(\beta,d)$ is a constant in the exponent of the Nash-Aronson estimate for parabolic Green's function (see Lemma 4.5 in \cite{Abdulle_Arjmand_Paganoni_18b} for details), which depends only on $d$ and $\beta$}. In addition, the choices
\begin{equation*}
L = k_{o} R, \quad T  = k_T R, 
\end{equation*}
with $0 < k_{o} < 1$, and $k_T = \sqrt{{\frac{c_2}{c_1}}} (1-k_o)$ results in the following convergence rate in terms of $R$
\begin{equation*}
\|  a^{0}_{T,R,L} -  a^{0} \|_{F}  \leq C \left( R^{-q-\frac{1}{2}} +\ncorr{R^{d/2}} e^{-\sqrt{c_1c_2} (1-k_{o})  R} \right),
\end{equation*}
where $C = \max\{C_1,C_2,C_3 \}$.  
\end{theorem}
\begin{Remark}
\ncorr{From the last result of \cref{Thm_Main_Result}, the error estimate can be decomposed into two terms, one exponentially decaying with respect to $R$ and the other one decaying as $R^{-q-\frac{1}{2}}$, where $q$ can be arbitrarily chosen. When $q$ has a finite value and $R$ is sufficiently large, the algebraic component of the error is dominating. However, it is possible to consider infinite-order filters (e.g. $\mu(x) = \frac{e^\phi(x)}{\int_{-1(/2)}^{1/2}e^{\phi(x)}\,dx} \mathds{1}_{[-1/2, 1/2]}(x)$, with $\phi(x)=\frac{2}{1-2|x|}$), or $q$ can be chosen as an increasing function of $R$. Both these strategies provide a decay rate faster than any algebraic convergence rate, hence an \emph{exponential} convergence of the error.}
\end{Remark}
In \cref{Thm_Main_Result}, The error $\sqrt{T} L^{-(q+1)}$ is the averaging error, which is obtained by using a filter $\mu_L\in \mathbb{K}^{q}(K_L)$. 
\corr{The order $q$ of the filter can be chosen arbitrarily large with no additional computational cost. This allows to have better convergence rates for the resonance error. However, for higher order filters we witness a \textit{plateau} in the convergence plot of the error, which is not present for low order filters, e.g., see \cref{Fig_TwoD_Layered}.}
The error $e^{-c_1 T}$ is related to the solution of the parabolic PDE \eqref{Eqn_Parabolic_CellProblem} for a finite $T$. Note that the parabolic PDE \eqref{Eqn_Parabolic_CellProblem} is introduced only for the analysis, but in practice, we don't solve it. The term $e^{-c_2 \frac{\left| R-L \right|^2}{T}}$ along with its prefactor is an upper bound for the boundary error, and it will decay exponentially fast  with respect to $R$ only if $T < |R-L|^2$.  

\begin{proof} We prove this theorem in several steps:

{\bf Step 1. Error decomposition.} 
The aim here is to show that the error can be split as
\begin{equation} \label{Ineq_Error_Splitting}
\| a^{0}_{T,R,L} - a^{0} \|_{F} \leq \mathcal{E}_{\text{av}} + \mathcal{E}_{\text{boundary}} + \mathcal{E}_{\text{truncation}}.
\end{equation}
The term $\mathcal{E}_{\text{av}}$ is the averaging error which decreases by using filters $\mu_L \in \mathbb{K}^{q}(K_L)$ with higher values for $q$. The error $\mathcal{E}_{\text{truncation}}$ is associated with truncation in time of the solutions of parabolic cell-problems. The boundary error $\mathcal{E}_{\text{boundary}}$ quantifies the effect of boundary conditions. 
To see this, we use \cref{Thm_Time_Integration} and write
\begin{align*}
e_i \cdot a^{0}_{T,R,L} e_j &:= \int_{K_L}  a_{ij}(\bx) \mu_L(\bx) \; d\bx + \int_{K_L} a_{ik}(\bx) \partial_{x_k} \chi_{T,R}^{j}(\bx) \mu_L(\bx)  \; d\bx \\
&=  \int_{K_L} a_{ij}(\bx)  \mu_L(\bx) \; d\bx + \int_{0}^{T} \int_{K_L}  a_{ik}(\bx) \partial_{x_k} u^{j}(t,\bx)  \mu_L(\bx) \; d\bx \;dt,
\end{align*}
where $u^{j}$ is the solution of the parabolic cell problem \eqref{Eqn_Parabolic_CellProblem}. In the same way, by \cref{Thm_Time_Integration_Periodic}, the exact homogenized coefficient given by \eqref{Eqn_Exact_Hom_Coeff} can be rewritten as 
\begin{align*}
e_i\cdot a^{0} e_j &= \int_{K} a_{ij}(\bx) \; d\bx + \int_{K} a_{ik}(\bx) \partial_{x_k} \chi^{j}(\bx) \; d\bx \\
&=  \int_{K} a_{ij}(\bx) \; d\bx + \int_{0}^{\infty} \int_{K} a_{ik}(\bx) \partial_{x_k} v^{j}(t,\bx) \; d\bx \;dt,
\end{align*}
where $v^{j}$ is the periodic parabolic solution in \eqref{Eqn_Parabolic_Periodic_CellProblem}, and $\chi^{j}$ is the solution to the periodic cell problem \eqref{Eqn_Periodic_Elliptic_CellProblem}. We exploit this equality to further decompose the error $\mathcal{E}:=\left| e_i (a^{0}_{T,R,L} -a^{0}) e_j \right|$ as follows
\begin{multline}\label{Eqn_Error_Splitting_Explicit_Terms}
\left| e_i \cdot \left( a^{0}_{T,R,L} -  a^{0} \right) e_j \right|  
\le  \underbrace{\left|
	\int_{K_L}  a_{ij}(\bx) \mu_L(\bx)\; d\bx 
	- \int_{K} a_{ij}(\bx) \; d\bx 
	\right|}_{\mathcal{E}_{\text{av1}}}  \\ 
+ \underbrace{\left| 
	  \int_{0}^{T} \int_{K_L} a_{ik}(\bx) \partial_{x_k} u^{j}(t,\bx) \mu_L(\bx)\; d\bx \;dt 
	- \int_{0}^{T} \int_{K_L} a_{ik}(\bx) \partial_{x_k} v^{j}(t,\bx) \mu_L(\bx) \; d\bx \;dt 
	\right|}_{\mathcal{E}_{\text{boundary}}}  \\ 
+ \underbrace{\left| 
      \int_{0}^{T} \int_{K_L}  a_{ik}(\bx) \partial_{x_k} v^{j}(t,\bx) \mu_L(\bx)\; d\bx \; dt 
	- \int_{0}^{T} \int_{K} a_{ik}(\bx) \partial_{x_k} v^{j}(t,\bx) \; d\bx \;dt 
	\right|}_{\mathcal{E}_{\text{av2}}}  \\ 
+ \underbrace{\left| 
	  \int_{0}^{T} \int_{K}  a_{ik}(\bx) \partial_{x_k} v^{j}(t,\bx) \; d\bx \; dt  
	- \int_{0}^{\infty} \int_{K}  a_{ik}(\bx) \partial_{x_k} v^{j}(t,\bx) \; d\bx \; dt 
	\right|}_{\mathcal{E}_{\text{truncation}}} 
\end{multline}
The averaging error in the splitting \eqref{Ineq_Error_Splitting} is then defined as $\mathcal{E}_{\text{av}} :=  \mathcal{E}_{\text{av1}} + \mathcal{E}_{\text{av2}}$. In the following steps we give bounds for all the errors. 

{\bf Step 2. The bound for $\mathcal{E}_{av}$.} The main result in this step is summarised in the following lemma. 
\begin{lemma}\label{Lem_E_av} 
Let $a \in \mathcal{M}_{per}(\alpha,\beta,K)$, $a e_i \in H_{\text{div}}(K)$ and $\mathcal{E}_{av}:= \mathcal{E}_{av1} + \mathcal{E}_{av2}$, where $\mathcal{E}_{av1}$ and $\mathcal{E}_{av2}$ are defined in \eqref{Eqn_Error_Splitting_Explicit_Terms}. Then 
\ncorr{\begin{align*}
\mathcal{E}_{av} \leq 
\begin{cases} C_1(\beta,\alpha,\mu) \sqrt{T}L^{-q-1} \| \nabla \cdot a e_j \|_{L^2(K)}& \text{ if }  \nabla \cdot a e_j \in L^2(K), \\
 C_2(\beta,\alpha,\mu) L^{-q-1}  \| \nabla \cdot a e_j \|_{W^1_{per}(K)} & \text{ if } \nabla \cdot a e_j \in W^1_{per}(K), 
\end{cases} 
\end{align*}
where $C_1 = C(\mu^{(q)}) \frac{\beta}{\alpha}$, and $C_2 = C(\mu^{(q)}) (\frac{\beta}{\alpha})^{3/2}$, and $C(\mu^{(q))}$ is the Lipschitz constant of $q$th derivative of the filter $\mu$.}

\end{lemma} 
\begin{proof}
By \cref{Lem_Filtering}, we can immediately see that 
\begin{align*}
\mathcal{E}_{\text{av1}}&:= \left|\int_{K_L}  a_{ij}(\bx) \mu_L(\bx)  \; d\bx - \int_{K} a_{ij}(\bx) \; d\bx \right| \\&\leq C L^{-q-1} \| a_{ij} \|_{L^{2}(K)} \leq C \beta L^{-q-1},
\end{align*}
\ncorr{where $C = C(\mu^{(q)})$ is the Lipschitz constant of $q$-th derivative of the filter $\mu$.}
Moreover, 
\begin{align*}
\mathcal{E}_{\text{av2}} &:= \left| \int_{0}^{T} \left(\int_{K_L}  a_{ik}(\bx) \partial_{x_k} v^{j}(t,\bx) \mu_L(\bx)\; d\bx - \int_{K} a_{ik}(\bx) \partial_{x_k} v^{j}(t,\bx) \; d\bx \right)\;dt \right| \\
&\leq C(\mu^{(q)}) L^{-q-1}\int_{0}^{T} \| a_{ik} \partial_{x_k} v^{j}(t,\cdot) \|_{L^2(K)} \; dt \\ &\leq C(\mu^{(q)}) L^{-q-1} \beta \int_{0}^{T} \| \nabla v^{j}(t,\cdot) \|_{L^2(K)} \; dt. 
\end{align*}	
If the tensor $a(\bx)$ has higher regularity, i.e. $\nabla \cdot ae_j \in W^1_{per}(K)$, we can directly estimate $\| \nabla v^{j} \|_{L^1(0,T;L^2(K))} := \int_{0}^{T} \| \nabla v^{j}(t,\cdot) \|_{L^2(K)} \; dt$ by \cref{Eq2_Estimate_Periodic_Parabolic_Problem} in \cref{Lem_Estimate_Periodic_Parabolic_Problem} and obtain
\ncorr{
\begin{equation*}
\mathcal{E}_{\text{av2}} \leq 
C(\mu^{(q)}) (\dfrac{\beta}{\alpha})^{3/2} L^{-q-1} \| \nabla \cdot ae_j \|_{W^1_{per}(K)}.
\end{equation*} }
Otherwise, if $\nabla \cdot ae_j \in L^2(K)$ only, we will apply Cauchy-Schwarz inequality which yields $$\int_{0}^{T} \| \nabla v^{j}(t,\cdot) \|_{L^2(K)} \; dt \leq \sqrt{T} \| \nabla v^{j} \|_{L^2(0,T;L^2(K))}.$$ 
Then employing \cref{Eq1_Estimate_Periodic_Parabolic_Problem} in \cref{Lem_Estimate_Periodic_Parabolic_Problem}, we obtain 
\ncorr{
\begin{equation*}
\mathcal{E}_{\text{av2}} \leq 
C(\mu^{(q)}) \dfrac{\beta}{\alpha} \sqrt{T} L^{-q-1} \| \nabla \cdot ae_{j} \|_{L^2(K)}.
\end{equation*} }
This completes the proof of the Lemma.
\end{proof}

{\bf Step 3. The bound for $\mathcal{E}_{\text{truncation}}$.}
\begin{lemma} Let $a \in \mathcal{M}_{\text{per}}(\alpha,\beta,K)$ and $a e_i \in H_{\text{div}}(K)$. Then the truncation error $\mathcal{E}_{\text{truncation}}$ defined in \eqref{Eqn_Error_Splitting_Explicit_Terms} satisfies the estimate 
\ncorr{
\begin{align*}
\mathcal{E}_{\text{truncation}} \leq  C(d,\alpha) e^{- \frac{\alpha \pi^2}{d} T} \|\nabla \cdot a e_j \|^{2}_{L^2(K)},
\end{align*}
where $C(d,\alpha) = \frac{d}{\alpha \pi^2}$. }
\end{lemma}

\begin{proof}
By using integration by parts and the Cauchy-Schwarz inequality we have
\begin{align*}
\mathcal{E}_{\text{truncation}} &:= \left| \int_{T}^{\infty} \int_{K}  a_{ik}(\bx) \partial_{x_k} v^{j}(t,\bx) \; d\bx \; dt \right| \\
&= \left| \int_{T}^{\infty} \int_{K}  \left( \partial_{x_k}a_{ik}(\bx) \right)  v^{j}(t,\bx) \; d\bx \; dt \right| \\
&\leq \int_{T}^{\infty} \| \nabla \cdot ae_i \|_{L^2(K)} \| v^{j}(t,\cdot) \|_{L^2(K)} \; dt \\
&\leq \| \nabla \cdot ae_i \|_{L^2(K)} \int_{T}^{\infty} e^{-\lambda_0 t} \| \nabla \cdot ae_j \|_{L^2(K)} \; dt \\
&= \| \nabla \cdot ae_j \|_{L^2(K)} \| \nabla \cdot ae_i \|_{L^2(K)} \dfrac{1}{\lambda_0} e^{-\lambda_0 T},
\end{align*}
\ncorr{Similar to the proof of Lemma \ref{Lem_Estimate_Periodic_Parabolic_Problem}}, we complete the proof by observing that $\lambda_0 \geq \frac{\alpha\pi^2}{d}$.  
\end{proof}

{\bf Step 4. The bound for $\mathcal{E}_{\text{boundary}}$.}
\begin{lemma}
Let $a \in \mathcal{M}_{\text{per}}(\alpha,\beta,K) $ and $ae_j \in H_{div}(K_R)$ for any $j = 1,\dots,d$, ${\mu_L \in \mathbb{K}^{q}(K_L)}$ with $L<\tilde{R}$, where $\tilde{R}$ is the largest integer such that $\tilde{R}\le R-1/2$. Then, the boundary error $\mathcal{E}_{\text{boundary}}$ defined in \eqref{Eqn_Error_Splitting_Explicit_Terms} satisfies the estimate
\ncorr{
\begin{align*} 
\mathcal{E}_{\text{boundary}} \leq C\|\nabla \cdot a\mathbf{e}_i \|^{2}_{L^2(K)}  \left(1 + \frac{T}{|R-L|}\right) \frac{R^{d-1}}{T^{d/2-1}} e^{-c_2 \frac{\left| R-L \right|^2}{T}},
\end{align*}
where $C = \gamma \dfrac{\beta^2}{\alpha \sqrt{\alpha c}}$, for a $\gamma>0$ indepedent of $a,R,L,T$, and  $0< c_2 <\nu(\beta,d)$, where $\nu(\beta,d)$ is the exponent in the Nash-Aronson estimate for parabolic Green's function (see Lemma 4.5 in \cite{Abdulle_Arjmand_Paganoni_18b} for details), which depends only on $d$ and $\beta$.}
\end{lemma}

\begin{proof}
To estimate the boundary error, we define $\theta^{j} = u^{j} - \rho v^{j}$, where the smooth function $\rho \in C^{\infty}_{c}(K_{R})$ satisfies
\begin{equation*}
\rho(\bx)  = \begin{cases} 1, & \text{ if } \bx \in K_{\tilde{R}}, \\
0, & \text{ on } \partial K_{R}.
\end{cases}
\end{equation*}
Then, $\left(u^{j} - v^{j} \right)(t,\bx)= \theta^{j}(t,\bx)$ for any $t>0$ and $\bx\in K_L\subset K_{\tilde{R}}$, hence
\begin{multline*}
\mathcal{E}_{\text{boundary}}  := \left| \int_{0}^{T} \int_{K_L}  a_{ik}(\bx) \partial_{x_k} \left( u^{j} - v^{j} \right)(t,\bx) \mu_L(\bx)\; d\bx \; dt \right|\\ 
= \left| \int_{0}^{T} \int_{K_L}  a_{ik}(\bx) \partial_{x_k} \theta^{j}(t,\bx) \mu_L(\bx) \; d\bx \;dt \right| \\
= \left| \int_{0}^{T} \int_{K_L} \partial_{x_k}\left( a_{ik}(\bx) \mu_L(\bx) \right)  \theta^{j}(t,\bx)  \; d\bx \; dt \right|.
\end{multline*}
Next, it follows that 
\begin{align*}
\mathcal{E}_{\text{boundary}} &\leq   \int_{K_L} | \partial_{x_k} \left( a_{ik}(\bx) \mu_L(\bx) \right)| \; d\bx \sup_{\bx \in K_{L}} \int_{0}^{T} | \theta^{j}(t,\bx) | \; dt \\
& \leq \| \mu \|_{W^{1,2}(K_L)} \| a e_{j} \|_{H^1_{div}(K_L)} \sup_{\bx \in K_{L}} \int_{0}^{T} | \theta^{j}(t,\bx) | \; dt\\
& \leq  C_{\mu} L^{-d/2} L^{d/2} \| a e_{j} \|_{H^1_{div}(K)} \sup_{\bx \in K_{L}} \int_{0}^{T} | \theta^{j}(t,\bx) | \; dt.
\end{align*}
Moreover, we use the following lemma from \cite{Abdulle_Arjmand_Paganoni_18b}, which gives a pointwise estimate for the function $\theta^{j}$. 
\begin{lemma} \ncorr{Let $a \in \mathcal{M}_{\text{per}}(\alpha,\beta,K) $, $ae_j \in W^{1,\infty}(K_R)$ for $j = 1,\dots,d$, and $0<L<R-2$. Let $\theta^j := u^{j} - \rho v^{j}$, where $u^{j}\in L^2\left( (0,T);H^1_0(K_R)\right)$ and $v^{j}\in L^p\left((0,T);W^{1,p}_{per}(K)\right)$, for $p>\max\lbrace\frac{d+2}{2},2\rbrace$ are the solutions of \eqref{Eqn_Parabolic_CellProblem} and \eqref{Eqn_Parabolic_Periodic_CellProblem} respectively. Then, there exist constants $c_2,C_1,C_2>0$, independent of $R$ and $L$, such that, for any $0<t<2c_2|R-L|^2$,
	\begin{equation*}
		\|\theta^i(\cdot,t)\|_{L^{\infty}(K_L)} \le 
		C_1\|\nabla \cdot\left(a(\cdot)\mathbf{e}_i\right)\|_{L^2(K)}
		\left( 1 + C_2 \frac{t}{|R-L|} \right)
		\frac{R^{d-1}}{t^{d/2}}
		e^{-c_2\frac{|R-L|^2}{t}}.
	\end{equation*} 
 Here  $C_1$ is a constant independent of $a,R,L,T$, and $C_2 = C \dfrac{\beta^2}{\alpha \sqrt{\alpha c_2}}$, for a $C>0$ indepedent of $a,R,L,T$. Moreover, the constant $0< c_2 <\nu(\beta,d)$, where $\nu(\beta,d)$ is the exponent in the Nash-Aronson estimate for parabolic Green's functions, and depends only on $d$ and $\beta$; see the proof of Lemma $4.5$ in \cite{Abdulle_Arjmand_Paganoni_18b} for further details.
}
\end{lemma}
\ncorr{Invoking the $L^1-L^{\infty}$ H\"older inequality and the monotone growth in $[0, \frac{2c_2|R-L|^2}{d}]$ of the function $f:t\mapsto\frac{1}{t^{d/2}} e^{-c_2\frac{|R-L|^2}{t}}$ we obtain}
\ncorr{
\begin{align*}
	\sup_{\bx \in K_L}& \int_{0}^{T}  \left| \theta^{j}(t,\bx) \right| \; dt  \\ 
	&\leq C_1\|\nabla \cdot\left(a(\cdot)\mathbf{e}_i\right)\|_{L^2(K)} R^{d-1}
	 \int_{0}^{T} \left( 1 + C_2 \frac{t}{|R-L|} \right)
	 \frac{1}{t^{d/2}}
	 e^{-c_2\frac{|R-L|^2}{t}} \; dt \\
	&\leq C_1\|\nabla \cdot\left(a(\cdot)\mathbf{e}_i\right)\|_{L^2(K)} R^{d-1}
	\int_{0}^{T} \left( 1 + C_2 \frac{t}{|R-L|} \right) \; dt
	\max_{t\in[0,T]}
	\frac{1}{t^{d/2}}
	e^{-c_2\frac{|R-L|^2}{t}} \\
	&\leq C_1\|\nabla \cdot\left(a(\cdot)\mathbf{e}_i\right)\|_{L^2(K)} R^{d-1}
	 \left( 1 + C_2 \frac{T}{2|R-L|} \right) 
	\frac{1}{T^{d/2-1}}
	e^{-c_2\frac{|R-L|^2}{T}} . 
\end{align*}
}
\end{proof}
\begin{Remark}
	We emphasize here that one of the key arguments in proving an exponentially decaying error bound for $\mathcal{E}_{\text{boundary}}$ is the requirement that $L < R$, see \cite{Abdulle_Arjmand_Paganoni_18b}. 
\end{Remark}
Collecting the results from Step 1 to Step 4 gives the bound of \cref{Thm_Main_Result}. 
\end{proof}
\section{Approximation of the exponential operator $e^{-TA}$}
\label{Sec_Arnoldi}
The exponential correction term $e^{-T A} g^j$ in the model problem \eqref{Eqn_New_CellProblem} needs to be approximated in computations. The very first approach would be to regard $e^{-T A} g^j$ as the solution (at time $T$) of a parabolic PDE with initial data $g^j$. Such a consideration would not lead to any gain in computational cost in comparison to the parabolic approach described in \cite{Abdulle_Arjmand_Paganoni_18b}. Here we describe two more efficient ways, based on \emph{spectral truncation} and a \emph{Krylov subspace iteration} (Arnoldi iteration), to compute the exponential correction term, which are far less expensive than solving a full parabolic PDE. We also show that the both approximations result in exponentially decaying errors bounds for increasing values of $R$, preserving the desired exponential decay in \cref{Thm_Main_Result}, and that the Arnoldi iteration is computationally less expensive than the spectral truncation.
\subsection{Spectral truncation}\label{SubSec_spect_trunc}
The correction term $e^{-AT} g^{j}$ in \eqref{Eqn_New_CellProblem} corresponds to the solution (at time $T$) of the parabolic PDE \cref{Eqn_Parabolic_CellProblem}.
A way of expressing the exponential operator is 
\begin{align*}\label{Eqn_CorrectionTerm_FullSum}
[e^{-AT} g^{j}](\bx):= \sum_{k=0}^{\infty} e^{-\lambda_k T} g^j_{k} \varphi_k(\bx),  \text{ where } \quad g^j_k := \langle g^j, \varphi_k \rangle_{L^2(K_{R})},
\end{align*}
\ncorr{where the operator $A$ is defined on $H^1_0(K_R)$. The eigenvalues and eigenfunctions of $A$ depend on domain $K_R$, but the $R$ subscript is omitted for brevity.}
If $T$ is not too small, most of the modes in the expansion can be neglected due to the exponential decay with respect to the eigenvalues. Hence solving a more expensive parabolic PDE can be avoided at the expense of computing a few dominant modes of the operator $A$. To this end, let  
\begin{equation*}
[e^{-A_N T} g^{j}](\bx):= \sum_{k=0}^{N-1} e^{-\lambda_k T} g^j_k \varphi_k(\bx).
\end{equation*}
Then the cell-problem \eqref{Eqn_New_CellProblem} can be approximated by 
\begin{equation} \label{Eqn_New_CellProblem_Truncated}
\begin{aligned}
-\nabla \cdot \left( a(\bx) \nabla \chi_{T,R,N}^{j}(\bx) \right)  &= g^{j}(\bx) - [e^{-A_{N} T} g^{j}](\bx)& \quad &\text{ in } K_{R} \\
\chi_{T,R,N}^{j}(\bx) &= 0 &\quad &\text{ on } \partial K_{R}.  
\end{aligned}
\end{equation}
Similarly, the homogenized coefficient \eqref{Eqn_New_Approx_Hom_Coeff} is approximated by
\begin{align} \label{Eqn_New_Approx_Truncated_Hom_Coeff}
e_i \cdot a^{0}_{T,R,L,N} e_j = \int_{K_{L}}  \Big( a_{ij}(\bx ) + a_{ik}(\bx) \partial_{x_k} \chi_{T,R,N}^{j}(\bx) \Big) \mu_L(\bx) \; d\bx. 
\end{align}
The spectrally truncated cell-problem \eqref{Eqn_New_CellProblem_Truncated} and the homogenized coefficient \eqref{Eqn_New_Approx_Truncated_Hom_Coeff} are the ultimate approximations used in computations. In the following lemma, we give a bound for the difference between $a^{0}_{T,R,L}$ and $a^{0}_{T,R,L,N}$ defined in \eqref{Eqn_New_Approx_Hom_Coeff} and \eqref{Eqn_New_Approx_Truncated_Hom_Coeff} respectively.

\begin{lemma} \label{Lem_Ahom_Diff_Trunc_Estimate} Let $a \in \mathcal{M}(\alpha,\beta,K_{R}),  a e_j \in H_{div}(K_{R})$, and $\mu_L \in \mathbb{K}^{q}(K_L)$. Moreover, let $a^{0}_{T,R,L}$ and $a^{0}_{T,R,L,N}$ be defined as in \eqref{Eqn_New_Approx_Hom_Coeff} and \eqref{Eqn_New_Approx_Truncated_Hom_Coeff} respectively. Then  
\begin{equation} \label{Estimate_Mode_Truncation}
\mathcal{E}_{\text{spectral}}:= | e_i \cdot \left( a^{0}_{T,R,L} - a^{0}_{T,R,L,N} \right) e_j| \leq  
C \left(\frac{R}{L}\right)^{\frac{d}{2}} R e^{-\frac{c_d N^{2/d} T}{R^2}}
\end{equation}
where $C(\alpha,\beta,d,\mu_L)$ and $c_d$ are constants independent of $T,R,L,N$.
\end{lemma}
\begin{proof}
Let 
\begin{align*}
[e^{-AT} g^{j}](\bx)  = \sum_{k=0}^{\infty} e^{-\lambda_k T} g^{j}_k \varphi_k(\bx), \quad   [e^{-A_{N}T} g^{j}](\bx)  = \sum_{k=0}^{N-1} e^{-\lambda_k T} g^{j}_{k} \varphi_k(\bx),
\end{align*}
where $\{ \lambda_j, \varphi_j(\bx) \}_{j=0}^{\infty}$ are the eigenvalue-function pairs of the operator $A = - \nabla \cdot \left( a  \nabla \right)$ with Dirichlet boundary conditions on the domain $K_R$. Moreover, let $E_{N}(\bx) :=  [e^{-AT} g^{j} - e^{-A_NT} g^{j}](\bx)$, with $g^{j}:=\nabla \cdot a e_j$. The eigenvalues of second order symmetric elliptic operators satisfy 
\begin{equation}\label{Ineq_Lambdak}
\lambda_k \geq c_d k^{2/d} |K_R|^{-2/d} = c_d k^{2/d} R^{-2},
\end{equation} 
where $c_d$ is a constant that depends on the dimension\footnote{The constant $c_d$ may depend on $\alpha$ and $\beta$ too. The value of $c_d$ can be approximated by computing a few eigenvalues $\lambda_k$ and finding the largest constant so that the relation \eqref{Ineq_Lambdak} holds.} $d$ and the ellipticity constant $\alpha$, see \cite{Li_Yau_1983,Safarov_Vassiliev_1997}. Then 
%
\begin{align*}
\| E_{N} \|_{L^{2}(K_{R})}^{2}  &\leq \sum_{\ell,k=N}^{\infty} e^{-\frac{c_d (\ell^{2/d}+k^{2/d}) T }{R^2}} g^{j}_{\ell} g_k^{j} \int_{K_{R}} \varphi_{\ell}(\bx) \varphi_k(\bx) \; d\bx  \\
&= \sum_{k=N}^{\infty} e^{-\frac{2 c_d k^{2/d} \,T }{R^2}} |g_k^{j}|^{2} 
\leq e^{-\frac{2 c_d N^{2/d} \,T }{R^2}} \| g^{j} \|^2_{L^2(K_{R})}.
\end{align*}
Taking the square root of both sides, we arrive at 
\begin{align*}
\| E_{N} \|_{L^2(K_{R})} \leq e^{- \frac{c_d N^{2/d} T}{R^{2}}}  \| g^{j} \|_{L^2(K_{R})}.
\end{align*}
Moreover, since the difference $\psi  := \chi^{j}_{T,R} -  \chi^{j}_{T,R,N}$ satisfies $-\nabla \cdot a(\bx) \nabla \psi(\bx) = E_N(\bx)$ with homogeneous Dirichlet BCs, standard elliptic regularity yields 
\begin{align*}
\| \chi^{j}_{T,R} -  \chi^{j}_{T,R,N} \|_{H^{1}_{0}(K_{R})} &\leq \frac{C_p(K_R)}{\alpha} \| E_N \|_{L^2(K_{R})} \\ &\leq C R e^{- \frac{c_d N^{2/d} T}{R^{2}}}  \| g^{j} \|_{L^2(K_{R})}\\
&\leq C R^{1+\frac{d}{2}} e^{- \frac{c_d N^{2/d} T}{R^{2}}} \| ae_j \|_{H_{div}(K)} ,
\end{align*}
where we have used the fact that the Poincar\'e constant $C_p(K_R)$ is bounded by $C_p(K_R) \leq  diam(K_R)/\pi=R 2^{1/d}/\pi$, see \cite{Payne_Weinerberger_1960}, and $\| g^{j} \|_{L^2(K_{R})} \leq |K_R|^{1/2} \| ae_j \|_{H_{div}(K)}$. Finally,
\begin{align*}
\left|  e_i \left( a^{0}_{T,R,L}  - a^{0}_{T,R,L,N} \right) e_j\right| &= \int_{K_L}  a_{ik}(\bx) \partial_{x_k} \left( \chi_{T,R} - \chi_{T,R,N} \right)(\bx) \mu_L(\bx) \; d\bx \\ 
& \ncorr{\leq \beta \|  \nabla \chi^{j}_{T,R} - \nabla \chi^{j}_{T,R,N} \|_{L^2(K_{R})} \| \mu_L \|_{L^{2}(K_R)} } \\
&\leq \beta   \|  \nabla \chi^{j}_{T,R} - \nabla \chi^{j}_{T,R,N} \|_{L^2(K_{R})} \frac{|K_{L}|^{1/2}}{L^d}\| \mu \|^{d}_{L^{\infty}(K)}\\ 
&\leq 
C \dfrac{R^{1+\frac{d}{2}}}{L^{\frac{d}{2}}} e^{-\frac{c_d N^{2/d} T}{R^2}}.  
\end{align*}
This completes the proof. 
\end{proof}

In \cref{Thm_Main_Result}, the optimal value for the parameter $T$ is $T = O(R)$. In order to get an exponential decay rate, such as $e^{-c R}$ for some positive $c$, in \cref{Lem_Ahom_Diff_Trunc_Estimate}, we then need to compute $N = O(R^d)$ eigenmodes. This growth of the number of eigenmodes with respect to the dimension is the main drawback of the naive spectral truncation leading to a high computational burden in higher dimensions. Therefore, in the next subsection we propose a much more efficient method based on the Krylov subspace iteration, and we show that the cost of the method will scale linearly in terms of the number of degrees of the freedom, while retaining the desired exponential accuracy for the approximation of the homogenized coefficient.

\begin{Remark}\label{Rem_Per_Assumption_Spectral}
Note that \cref{Lem_Ahom_Diff_Trunc_Estimate} does not assume the periodicity of $a$, since the proof only relies on the decay of the eigenvalues of general second order elliptic operators. 
\end{Remark}

\subsection{Approximation by the Arnoldi method}\label{SubSec_Arnoldi}
In order to introduce the approximation by the Arnoldi method, we consider a discretization $\mathcal{A}_{h}$ of the operator $-\nabla \cdot \left( a \nabla \right)$ e.g., by a second order centred finite difference scheme, where $h$ is a discretization parameter. Moreover, assume that the size of the matrix $\mathcal{A}_h$ is $N\times N$, and that $g_h \in \mathbb{R}^{N}$ is a finite dimensional representation of $g^j$ on a uniform computational grid\footnote{Note that this is not the same parameter $N$ as in \cref{SubSec_spect_trunc}}. For the sake of simplicity,  in this section we will ignore the $j$ superscript. Denoting $F(z) = e^{-zT}$, the idea behind the Arnoldi algorithm is to look for an approximation for $F(\mathcal{A}_h) g_{h}$ starting by a unitary transformation of $\mathcal{A}_h$ in the form $\mathcal{H}  = \mathcal{Q}^{*} \mathcal{A}_h \mathcal{Q}$, where $\mathcal{Q} \in \mathbb{R}^{N\times k}$, $\mathcal{H} \in \mathbb{R}^{k\times k}$ is an upper-Hessenberg matrix, and $k \ll N$ so that the matrix $\mathcal{A}_h$ is projected into a lower dimensional space. The term $F(\mathcal{A}_h) g_h$ can then be approximated by 
\begin{equation*}
F(\mathcal{A}_h) g_h \approx \mathcal{Q} F(\mathcal{H}) \mathcal{Q}^{*} g_h.
\end{equation*}
Therefore, computing the computationally expensive exponential matrix function $F(\mathcal{A}_h)$ of size $N\times N$ is avoided by instead computing $F(\mathcal{H})$, with a smaller computational cost. An important question that arises is in relation with the approximation error coming from the Arnoldi algorithm. The following theorem from \cite{Hochbruck_Lubich_97} provides an upper bound for such an approximation. 
\begin{theorem} (Hochbruck, Lubich \cite{Hochbruck_Lubich_97}) \label{Hochbruck_Lubich_Thm}
Let $\mathcal{B} \in \mathbb{R}^{N\times N}$ be a Hermitian negative semi-definite matrix with eigenvalues in $[-\rho,0]$ and set $\beta = \| g \|_2$ where $g\in \mathbb{R}^{N}$. Moreover, let $\mathcal{H} =\mathcal{Q}^{*} \mathcal{B} \mathcal{Q}$ be a unitary transformation of $\mathcal{B}$ via an Arnoldi procedure with $\mathcal{H} \in \mathbb{R}^{k\times k}$ and $\mathcal{Q} \in \mathbb{R}^{N\times k}$. Then the following estimate holds
\begin{equation}
\| e^{\mathcal{B}} g - \mathcal{Q} e^{\mathcal{H}} \mathcal{Q}^{*} g\|_2 \leq \begin{cases} 10 \beta e^{-4 k^2/(5\rho)}, & \sqrt{\rho} \leq k \leq \rho/2 \\
\dfrac{40 \beta}{\rho} e^{-\rho/4} \left( \dfrac{e \rho}{4k} \right)^k, & k \geq \rho/2.
\end{cases}
\end{equation}
\end{theorem}
\ncorr{%
\begin{lemma}\label{Lemma_spectralRadius}
	Let $\mathcal{A}_{h} \in \mathbb{R}^{N\times N}$ be a second order centred finite difference approximation, \ncorr{with step-size $h$}, of the operator $-\nabla \cdot \left( a \nabla \right)$ on $K_R$, and let $\rho\left(\mathcal{A}_h\right)$ denote its spectral radius. 
	Then, $\rho\left(\mathcal{A}_h\right)=O(h^{-2})$ and it is independent of $R$.
\end{lemma}
\begin{proof}
	Let us define $\hat{a}(\by) = a(\bx)$ and $\hat{u}(\by) = u(\bx)$, where $\by =\bx/R$.
	Then, 
	\begin{equation}\label{eq:scaling_eq}
		-\nabla_{\bx} \cdot \left( a(\bx) \nabla_{\bx} u(\bx) \right)  = - \frac{ 1}{R^2} \nabla_{\by} \cdot \left( \hat{a}(\by) \nabla_{\by} \hat{u}(\by) \right),
	\end{equation}
	where the operator $-\nabla \cdot \left( \hat{a} \nabla \right)$ is independent of $R$.
	Let $\hat{\mathcal{A}}_{\hat{h}} \in \mathbb{R}^{N\times N}$ be a second order centred finite difference approximation, with step-size $\hat{h}=h/R$, of $-\nabla \cdot \left( \hat{a} \nabla \right)$. Since \eqref{eq:scaling_eq} holds for any sufficiently smooth function $u$, it follows that 
	 \begin{equation*}
	 	{\mathcal{A}}_h = \frac{1}{R^2} \hat{\mathcal{A}}_{\hat{h}}.
	 \end{equation*}
 	It is known that the spectral radius of $\hat{\mathcal{A}}_{\hat{h}}$ scales as $\hat{h}^{-2}$ (see, e.g. \cite{LeVeque2007}, Chapter 3  for the proof in the case of the Laplacian operator), hence
 	$$
 	\rho(\mathcal{A}_h) =  \frac{1}{R^2} \rho(\hat{\mathcal{A}}_{\hat{h}}) = \frac{1}{R^2} O(\hat{h}^{-2}) = O(h^{-2}).
 	$$
\phantom{some fake test}
\end{proof}
}
\begin{corollary}\label{Cor_Arnoldi}
Let $\mathcal{A}_{h} \in \mathbb{R}^{N\times N}$, with $h\leq 1$, be a second order centred finite difference approximation of the operator  $-\nabla \cdot \left( a \nabla \right)$ and set $\beta = \| g_{h} \|_2$ where $g_{h}\in \mathbb{R}^{N}$. Moreover, let $\mathcal{H} =\mathcal{Q}^{*} \mathcal{B} \mathcal{Q}$ be a unitary transformation of $\mathcal{A}_{h}$ via an Arnoldi procedure with $\mathcal{H} \in \mathbb{R}^{k\times k}$, $\mathcal{Q} \in \mathbb{R}^{k\times N}$, and $F(z) = e^{- zT}$. Then the following estimate holds
\begin{equation} \label{Est_Arnoldi_2}
\| F(\mathcal{A}_h) g_h - \mathcal{Q} F(\mathcal{H}) \mathcal{Q}^{*} g_h\|_2 \leq 10 \beta e^{- 4 k^2/(5 c_d N^{2/d} R^{-2} T)},
\end{equation}
for $\sqrt{ c_d N^{2/d} R^{-2} T} \leq k \leq c_d N^{2/d} R^{-2} T/2$, where $c_d$ is a constant which depends only on the dimension. Moreover, when $k= \sqrt{c_d N^{2/d} R^{-2}} T/2$,\ncorr{and for sufficiently small $h$ (whenever $c_d < 1$)}, the estimate reads as
\begin{equation}\label{Est_Arnoldi}
\| F(\mathcal{A}_h) g_h - \mathcal{Q} F(\mathcal{H}) \mathcal{Q}^{*} g_h\|_2 \leq 10 \beta e^{- T/5},\quad \text{for }T\ge 4.
\end{equation}
\end{corollary}
\begin{proof}
The proof of \eqref{Est_Arnoldi_2} follows by a direct application of the \cref{Hochbruck_Lubich_Thm} and the fact that the spectral radius of the operator $\mathcal{A}_{h}T$ is given by (see \cref{Lemma_spectralRadius})
$$
\ncorr{\rho(\mathcal{A}_hT) = O(h^{-2} T)  = O(N^{2/d} R^{-2} T).}
$$
\ncorr{For the second part of the proof, it is sufficient to verify that the optimal choice $k= \sqrt{c_d N^{2/d} R^{-2}} T/2$ is within the interval  $\sqrt{ c_d N^{2/d} R^{-2} T} \leq k \leq c_d N^{2/d} R^{-2} T/2$. This is easy to see if we consider two cases: 1) If $c_d \geq 1$, then by using the fact that $h = R/N^{1/d} \leq 1$ and $T\ge 4$, the inequalities are satisfied. 2) When $c_d \leq 1$, then the inequalities are satisfied for sufficiently small $h$. }
\end{proof}
The advantage of using an approximation for the exponential correction term via the Arnoldi approach is that the number of basis functions required in the Arnoldi iteration is practically independent of the dimension of the problem. In other words, denoting the numbers of degrees of freedom in $d$-dimensions by  \ncorr{$N:=n^{d} \simeq (R/h)^d $, only $k= \sqrt{c_d} T/(2h)$} basis functions are needed to obtain an exponentially accurate approximation for the exponential correction $e^{-A T} g$ up to a discretization error, see the estimate \eqref{Est_Arnoldi}. \ncorr{An estimate for a fully discrete approximation of the homogenized coefficient can also be derived, following the lines of \cite{Abdulle_05}. In \cite{Abdulle_05}, the author identifies the so-called \emph{microscopic discretization error} as one of the error sources for FE-HMM. Such an error is due to the numerical approximation of the microscopic correctors $\chi^j_{T,R}$. If an $s$-th order FEM is used for the micro-problems, then the micro-error reads as $e_{MIC} = O( h^{2s})$. In our setting, apart from the micro error,  there are also resonance error and the error due to Arnoldi approximation.  We would then expect that the overall errors for a fully discrete analysis of approximation for the homogenized coefficient read as
\begin{equation*}
\mathcal{E} \leq C \left( h^{2s} + \mathcal{E}_{resonance} + \mathcal{E}_{Arnoldi} \right)
\end{equation*}
A fully discrete analysis is skipped in the present paper so as to remain in line with our main goal of proving an error bound for the approximation of the homogenized coefficient in a continuous setting. But we emphasize that it follows directly from \cite{Abdulle_05} and the results therein.}

\subsection{Approximation of the cell problem and computational cost} \label{SubSec_CompCost}
The Arnoldi iteration can be used in different ways to approximate the solution of the modified elliptic PDE \eqref{Eqn_New_CellProblem}. A standard finite element/difference discretization of the problem \eqref{Eqn_New_CellProblem} results in the following system\footnote{For simplicity all the indices are skipped in this discussion.} 
\begin{equation} \label{Eq_Main_Discrete}
\mathcal{A}_h  \chi_{h}= g_{h} - e^{-T \mathcal{A}_{h} } g_{h}.
\end{equation}
Here we present three different ways based on the Arnoldi iteration to solve \eqref{Eq_Main_Discrete}.

{\bf Approach 1.} Let $F_1(z)  = e^{-tz}$, then the system \eqref{Eq_Main_Discrete} can be approximated by
\begin{equation} \label{Eq_Discrete_Arnoldi_1}
\mathcal{A}_h  \tilde{\chi}_{h}= g_{h} - \mathcal{Q} F_1(\mathcal{H}) \mathcal{Q}^{*} g_{h}.
\end{equation}

{\bf Approach 2.} Let $F_2(z)  = 1- e^{-tz}$, then the system \eqref{Eq_Main_Discrete} can be approximated by
\begin{equation} \label{Eq_Discrete_Arnoldi_2}
\mathcal{A}_h  \tilde{\chi}_{h}= \mathcal{Q} F_2(\mathcal{H}) \mathcal{Q}^{*} g_{h}.
\end{equation}

{\bf Approach 3.} Let $F_3(z)  = z^{-1} (1- e^{-tz})$, then the system \eqref{Eq_Main_Discrete} can be approximated by
\begin{equation} \label{Eq_Discrete_Arnoldi_3}
\tilde{\chi}_{h}= \mathcal{Q} F_3(\mathcal{H}) \mathcal{Q}^{*} g_{h}.
\end{equation}

Assuming that the systems in approaches $1,2$ are inverted by a linearly scaling algorithm, such as the multigrid, the overall computational costs of all these three formulations are dominated by the Arnoldi iteration, where given the matrix $[\mathcal{A}_h]_{N\times N}$, the matrices $\mathcal{Q}$ and $\mathcal{H}$ are computed. The Arnoldi algorithm consists of an outer loop for $j=1:k$, where in total $k \ll N$ sparse matrix vector multiplications of the form $\mathcal{A}_h g_h$ are needed. Moreover, there is an orthogonalisation process which occurs at an inner loop for $i=1:j$, where the essential cost is due to a vector-vector multiplication $g_h^T g_h$ of two dense $N\times 1$ vectors. The overall cost of the Arnoldi iteration, exploiting the inherent sparsity of $\mathcal{A}_h$ becomes
\begin{equation*}
Cost_{Arnoldi}  \approx \sum_{j=1}^{k} \left( C_d N + 4 \sum_{i=1}^{j} N \right) = O(N k^2). 
\end{equation*} 
If $k$ is fixed \textit{a-priori} instead of following the scaling of \cref{Cor_Arnoldi}, then  the cost of the algorithm will grow linearly with $N$. A more rigorous analysis can be done by using the analysis in \cref{SubSec_Arnoldi}, where the optimal value of $k$ for the approach $1$ has been presented. Following the result of \cref{Cor_Arnoldi}, we find that \ncorr{$k^2 = O(N^{2/d})$}. Hence the overall cost becomes 
\ncorr{\begin{equation*}
Cost_{Arnoldi} = O(N k^2) = O(N^{1+2/d}). 
\end{equation*}}
\corr{Using the relation $N = R^d h^{-d}$, where $h$ is a fixed mesh size, we can write the computational cost of the modified and the standard elliptic upscaling approaches as a function of $R$ and $h$, see \cref{Table_Cost}. The global errors, which are composed of the resonance and the discretization errors, are also reported in \cref{Table_Cost}. 
	The resonance error scales as $R^{-q-1/2} $ for the modified elliptic approach, see \cref{Thm_Main_Result}, while it decays as $R^{-1}$ for the standard elliptic case, see \cref{eq_conv_standard_modelling_error}. The discretization error is assumed to be of order $O(h^s)$ in both cases.
	In order to derive the scaling of the cost with respect to the accuracy, we impose the global error to be smaller than a prescribed tolerance $tol$. So, for the modified elliptic case, we choose $R$ and $h$ such that $R^{-q-1/2}\approx tol$ and $h^{s}\approx tol$, while $R^{-1}\approx tol$ and $h^{s}\approx tol$ for the standard elliptic case.
}
Therefore, the modified elliptic approach has a lower cost to reach a certain tolerance $tol$ when 
\ncorr{\begin{equation*}
\dfrac{d+2}{q+1/2} + \dfrac{2}{s} < d,
\end{equation*}}
which is easily achieved by using filters with better regularity properties (large $q$), as well as high order numerical methods for the approximation of the elliptic PDE \eqref{Eqn_New_CellProblem}. 

Note that although an estimate for the difference between $\mathcal{Q} F_1(\mathcal{H}) \mathcal{Q}^{*} g_{h}$ and $e^{-T \mathcal{A}_{h}} g_{h}$ is available, see the results stated in \cref{SubSec_Arnoldi}, error estimates for more complicated matrix functions such as $F_3(z)  = z^{-1} (1- e^{-tz})$ used in approach $3$ above are not known, \cite{Higham_Book}. Nevertheless, from a computational point of view, the approach $3$ has a slight advantage of skipping the inversion of the large sparse matrix $\mathcal{A}_h$, and hence is used in the simulations of this paper. There is another approach based on the rational Krylov subspace iteration, see e.g., \cite{Guntel_13}, which better suits the treatment of functions such as $F_3$, but the method requires an inversion of $\mathcal{A}_h$ for each column of the matrix $Q$. 
\begin{table}[h]
	\centering
\begin{tabular}{ | c | c | c | c | }
	\hline
 Elliptic cell problem & Computational cost & Error & Cost($tol$) \\
 \hline
 Modified & $R^{2+d} h^{-d-2}$  & $R^{-q-1/2} + h^{s}$ & $tol^{-\frac{2d+4}{2q+1}-\frac{d+2}{s}}$   \\
 Standard & $R^{d} h^{-d}$ & $R^{-1} + h^s$ & $tol^{-d-\frac{d}{s}}$   \\
 \hline
\end{tabular}
\caption{Cost to reach a tolerance $tol$. The cost of solving the modified elliptic problem \eqref{Eqn_New_CellProblem} by Arnoldi approach 1 is compared to that of the standard elliptic problem \eqref{Eqn_Periodic_Elliptic_CellProblem}.} \label{Table_Cost}
\end{table} 
%
%
\section{Numerical tests}\label{Sec_Numerical_Results}
In this section, we provide examples in two dimensions to verify the theoretical results stated in \cref{Thm_Main_Result}. Moreover, additional numerical tests are provided to show that the method performs equally well even when the regularity and structural assumptions of the theorem are violated. In particular, the test cases include a periodic medium, a discontinuous layered medium, a quasi-periodic medium, as well as a random medium. These results are discussed in separate subsections below. \ncorr{Note that in all the simulations below, we will use the Arnoldi approach 3 for the numerical approximation of the modified elliptic approach, and we will use the Frobenius norm $\| a^{0} - a^{0}_{T,R,L}  \|_{F}$ to compute the errors.}
%

\begin{figure}[ht]%
    \centering
    \subfloat[$a(\bx) = diag(a_{11}(\bx),a_{22}(\bx))$]{
    	\includegraphics[height=4.65cm]{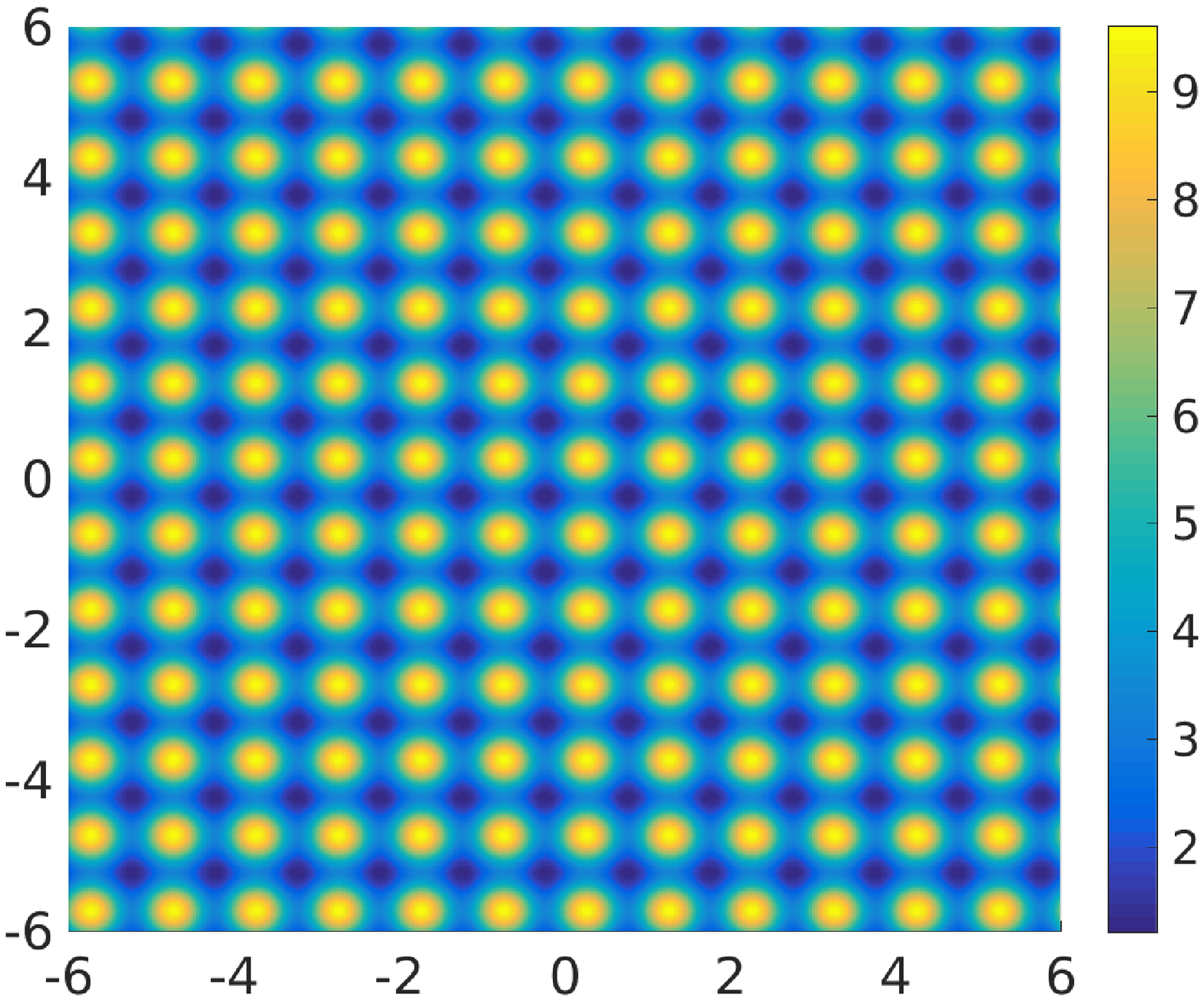} 
    }
    \subfloat[The upscaling error]{\includegraphics[height=4.65cm]{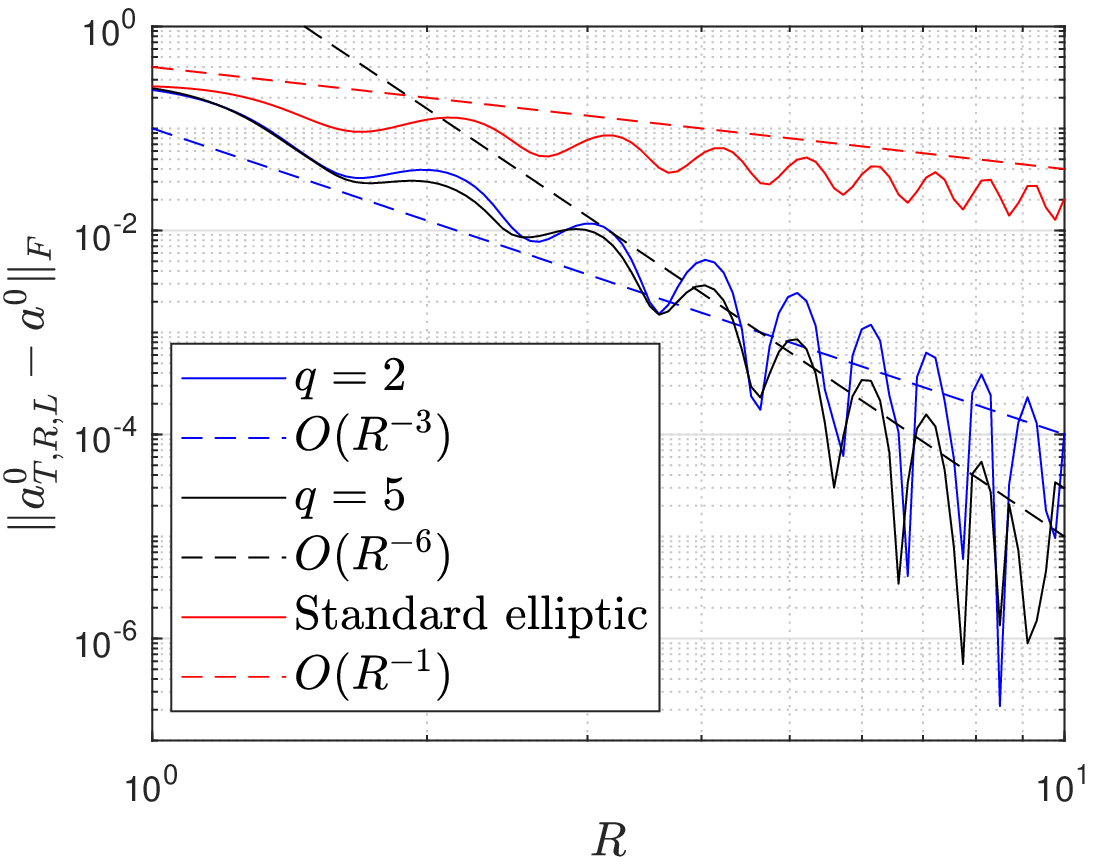}}%
    \caption{A two dimensional smooth medium \eqref{eqn_Coeff_Smooth2D} with $C_1=2.1, C_2 = 1$}%
    \label{Fig_TwoD_Ex1}%
\end{figure}

\textbf{Example 1. A smooth periodic coefficient.}
As our first example, we consider the following two-dimensional coefficient
\begin{align*}
a(\bx)  = \prod_{j=1}^2 \left( C_1 + C_2 \sin(2 \pi x_j)   \right) I,
\end{align*}
where $I$ is the $2\times 2$ identity matrix, see the left picture in \cref{Fig_TwoD_Ex1} for a graphical representation of $a$. \ncorr{Since the coefficients are separable, i.e. $a(\bx)=a_1(x_1)a_2(x_2)I$, the homogenized limit can be computed as (see \cite{Froese2009}, Section 4.2): 
\begin{align}\label{eqn_Coeff_Smooth2D}
a^{0} = \fint_K a_{11}(x) \,dx  \left( \fint_K a^{-1}_{22}(x) \,dx \right)^{-1} I = \left( C_1 \sqrt{ C_1^2  - C_2^2}  \right) I.
\end{align}}
In \cref{Fig_TwoD_Ex1}, the upscaling error $\|a^{0}_{T,R,L,N}  - a^{0}\|_{F}$ is shown for increasing values of $R$. The parameter values $T$ and $L$ are chosen optimally as stated in \cref{Thm_Main_Result}, with $k_0 = \frac{2}{3}$, $\alpha = \min_{\bx \in K} a(\bx)$,  $\beta = \max_{\bx \in K} a(\bx)$, and $k_T = \frac{\sqrt{d}}{2\pi\sqrt{\alpha\beta}}(1 - k_{o})$\ncorr{, since $c_1 = \frac{\alpha \pi^2}{d}$ (see \cref{Thm_Main_Result}) and $c_2$ can be approximated by $1/(4\beta)$ (see \cite{Abdulle_Arjmand_Paganoni_18b}, Section 5)}. The number of basis functions in Arnoldi algorithm to approximate the right hand side is $k=\min(700, N^{1/d})$ (where $N$ is the total number of degrees of freedom) for all values of $R$ since the Arnoldi's error is typically much smaller than the rest of the errors. Two different kernels with $q=2$ and $q=5$ are used in the simulations. The cell-problem \eqref{Eqn_New_CellProblem_Truncated} is approximated by a second order finite difference scheme with the stepsize $h= 1/120$. \ncorr{The total number of grid-points per direction is thus $R/h$}. The numerical results show that the overall error is dominated by the filtering error even for moderate values of $R$, and that arbitrarily high convergence rates are obtained by using kernels with better regularity properties. \ncorr{ Note that the filtering error can be estimated as $\mathcal{E}_{av} \leq C L^{-q-1}$ because $\nabla\cdot ae_j \in W^1_{per}(K)$ (see \cref{Lem_E_av}), which also explains the improved decay rate in the numerical results. Taking larger values for $q$ does not increase the overall computational cost, as the filters with larger $q$ values can be easily precomputed. Moreover, in the same figure, we also plot the corresponding error for the standard elliptic cell problem \eqref{Eqn_Approx_Periodic_Elliptic_CellProblem} with the same second order discretization with stepsize $h=1/120$. We observe a first order convergence rate for the resonance error when we use the standard elliptic problem.}

\begin{figure}[h]%
    \centering
    \subfloat[$a(\bx) = diag(a_{11}(\bx),a_{22}(\bx))$]{
    	\includegraphics[height=4.65cm]{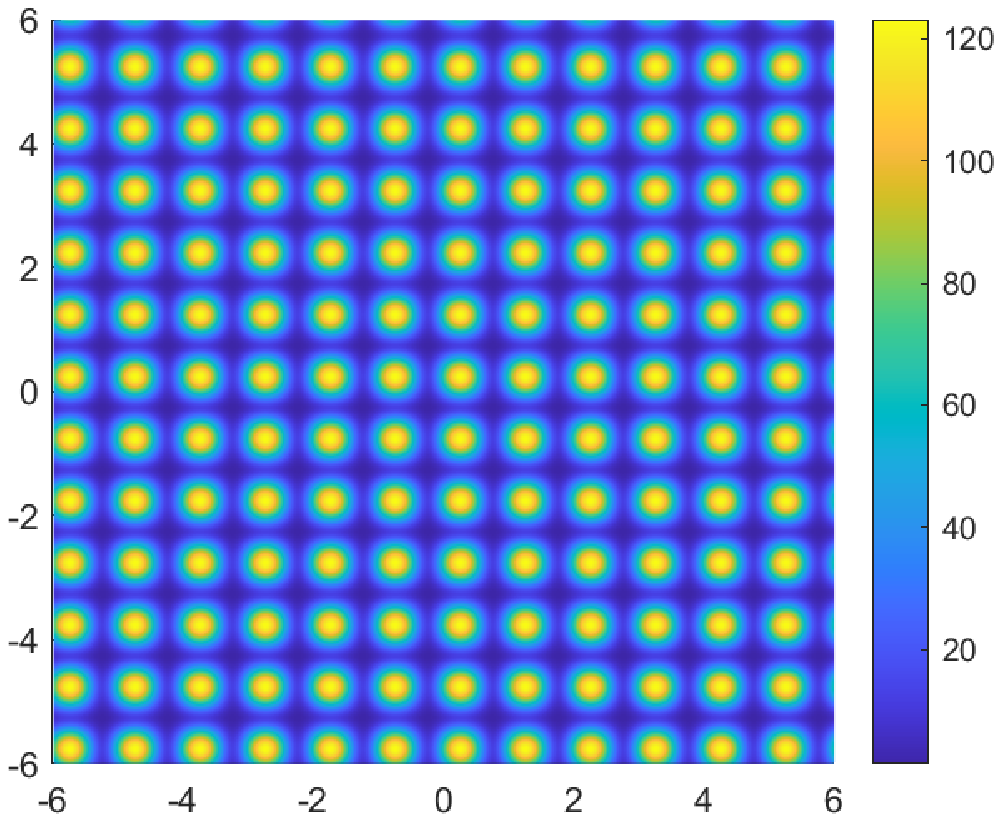} 
    }
    \subfloat[The upscaling error]{\includegraphics[height=4.65cm]{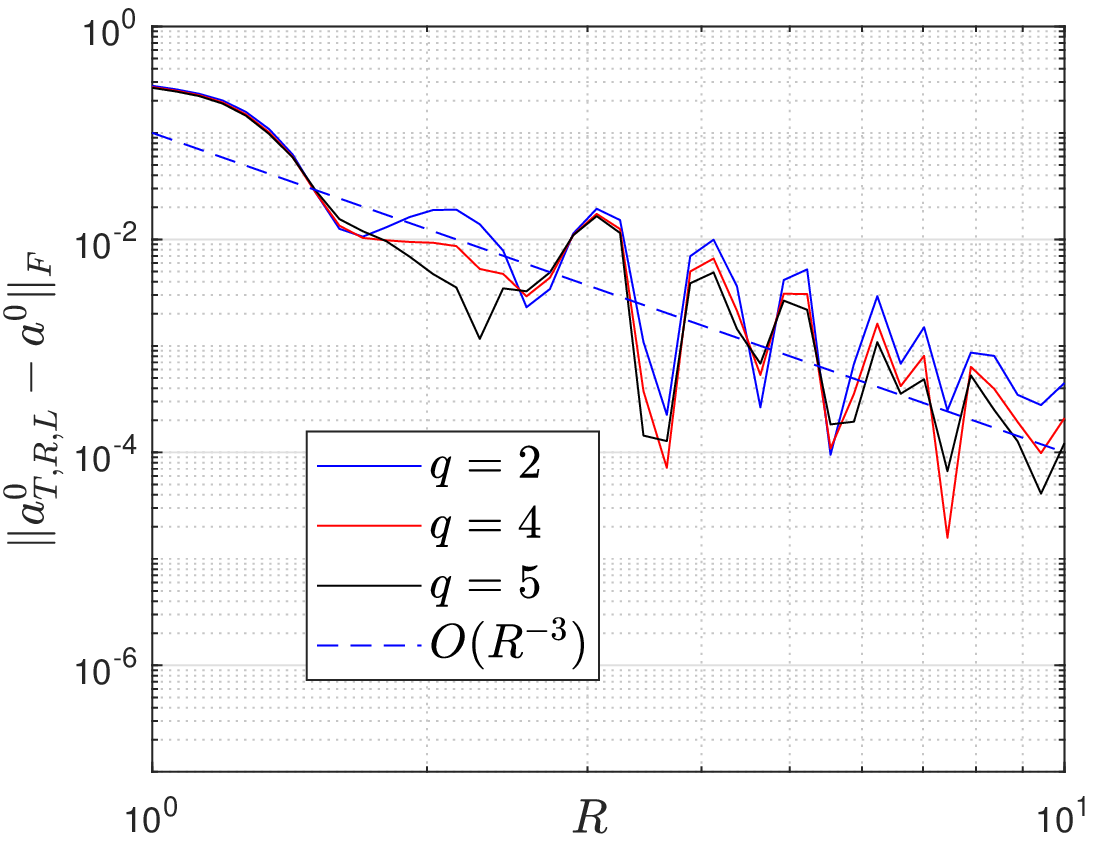}}%
    \caption{A two dimensional smooth medium \eqref{eqn_Coeff_Smooth2D} with contrast ratio $\approx 101$, $C_1=6.1, C_2 = 5$}%
    \label{Fig_TwoD_Smooth_HighContrast}%
\end{figure}

\ncorr{In Figure \ref{Fig_TwoD_Smooth_HighContrast}, we choose the constasts $C_1 = 6.1$ and $C_2 = 5$ to increase the contrast ratio of the coefficient to $\beta/\alpha \approx 101$. We use $h=1/50$, and choose optimal values for $k_T$ using the new values of $\alpha$ and $\beta$, and study the convergence with respect to increasing values of $R$. Other parameter values are chosen as in Figure \ref{Fig_TwoD_Ex1}. In this high contrast regime, we do not observe a strong dependency of the resonance error on the parameter $q$. This is expected since in the high contrast regime, the boundary error will dominate all other errors, and the effect of taking higher values for $q$ will be seen only for very large values of $R$; once the averaging error will be more dominant.}

\bigskip

%

\textbf{Example 2. A discontinuous periodic coefficient.} 
The second example is a layered medium characterised by the coefficient $a(\bx) = diag(a_{11}(x_1),a_{11}(x_1))$, where
\begin{align}\label{eqn_Coeff_Layered2D}
a_{11}(\bx)  = \begin{cases}  10  & 0\leq x_1 < \frac{1}{2} \\
1  & \frac{1}{2}\leq x_1 < 1.
\end{cases}
\end{align}
Such a choice is to test the generality of the method when the regularity assumption on the coefficient is relaxed. The exact homogenized coefficient is again constant and given by
\begin{align*}
a^{0} = diag(20/11, 11/2).
\end{align*}
All the numerical parameters are chosen identical to those in example $1$, with an obvious adaptation of $\alpha$ and $\beta$. Similar to example $1$, higher order convergence rates are achieved upon using higher order kernels, showing the generality of the method also for problems in discontinuous media.

%

\begin{figure}[t!]%
    \centering
    \subfloat[$a(\bx) = diag(a_{11}(\bx),a_{22}(\bx))$]{
%
%
\begin{tikzpicture}

\begin{axis}[%
width=5cm,
height=3.35cm,
scale only axis,
xmin=-6,
xmax=6,
xlabel style={font=\color{white!15!black}},
xlabel={$x\phantom{R}$},
xtick={-6,-4,-2 , 0,2,4,6},
xticklabels={-6,-4,-2 , 0,2,4,6},
ymin=-1.5,
ymax=9,
ylabel style={font=\color{white!15!black}},
axis background/.style={fill=white},
xmajorgrids,
ymajorgrids,
legend style={font=\scriptsize, at={(.02,.02)}, anchor = south west,legend cell align=left, align=left, draw=white!15!black}
]

\addplot[color=blue, domain=-6:6, samples=1024]
	plot (\x,{ 4 + cos(2*pi*\x r) + cos(2*pi*sqrt(2)*\x r) });
\addlegendentry{$a_{11}(0,x)$}

\addplot[color=black, domain=-6:6, samples=1024]
	plot (\x,{ 6 + (sin(2*pi*\x r))^2 + (sin(2*pi*sqrt(2)*\x r))^2 });
\addlegendentry{$a_{22}(x,x)$}

\end{axis}
\end{tikzpicture}%
}%
    \subfloat[The upscaling error]{\input{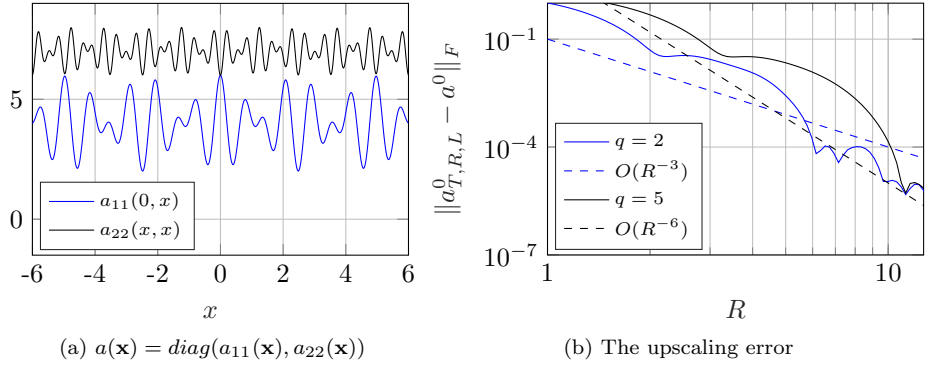} }%
    \caption{A two dimensional quasi-periodic medium \eqref{eqn_Coeff_QuasiPer2D}}%
    \label{Fig_TwoD_QP}%
\end{figure}

\textbf{Example 3. A quasi-periodic coefficient.} 
To test the applicability of the method beyond the periodic setting, we consider a quasi-periodic coefficient given by $a(\bx) = diag(a_{11},a_{22})$, where 
\begin{gather} \label{eqn_Coeff_QuasiPer2D}
a_{11}(\bx) = 4 + \cos(2 \pi (x_1 + x_2)) + \cos(2 \pi \sqrt{2} (x_1 + x_2)),\\
a_{22}(\bx) = 6 + \sin^2(2 \pi x_1 ) + \sin^2(2 \pi \sqrt{2} x_1). \nonumber
\end{gather}
The very same coefficient has been used also in the elliptic approach proposed in \cite{Gloria_2011}. In this paper, such a choice for the coefficient has been intentional  as it allows for a comparison between the two methods. In this particular setting, the homogenized coefficient is not easy to compute and therefore the value of $a^{0}_{T,R,L,N}$ with the largest $R$ is used instead of $a^{0}$ (similar to \cite{Gloria_2011}). All the parameter values are chosen identical to the example $1$ in this paper. \cref{Fig_TwoD_QP} shows a fast decay of the error down to $10^{-5}$ for moderate values of $R$, i.e., $R \approx 10$. It is worth mentioning that such an error tolerance is achieved only for $R \approx 20$ in the zero-order approach from \cite{Gloria_2011}.

%

\begin{figure}[t!]%
    \centering
       \subfloat[$a(\bx) = diag(a_{11}(\bx),a_{22}(\bx))$]{{\includegraphics[height=4.65cm]{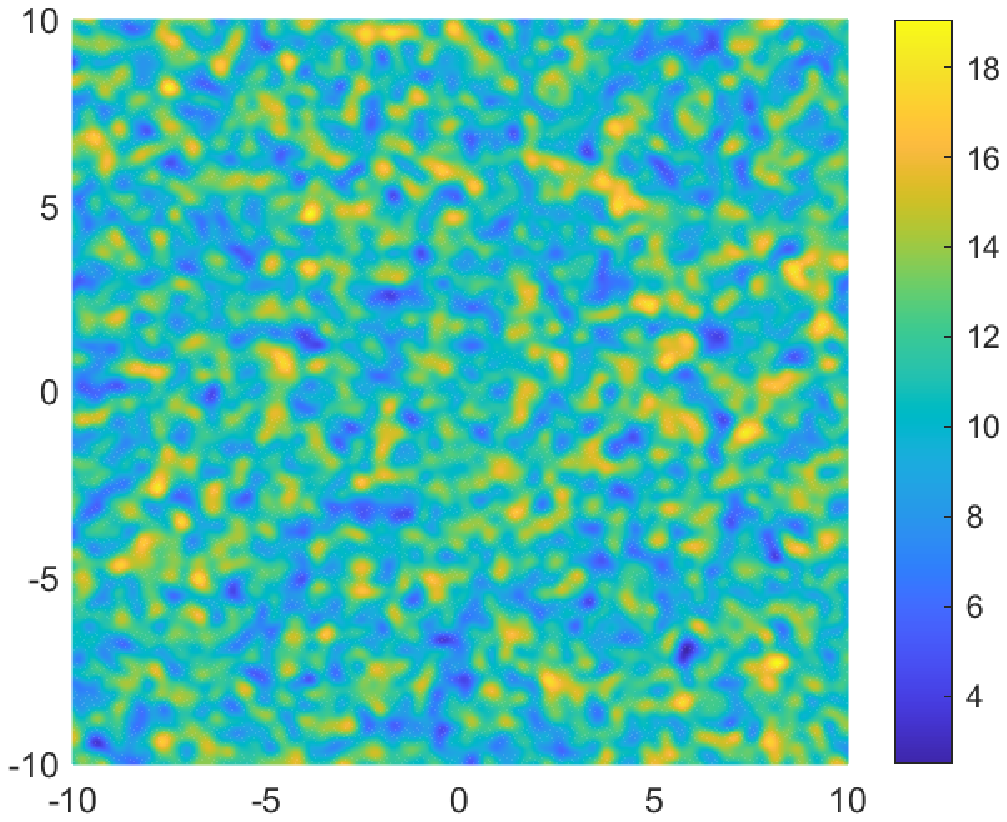} }}%
	\subfloat[The upscaling error]{{\includegraphics[height=4.65cm]{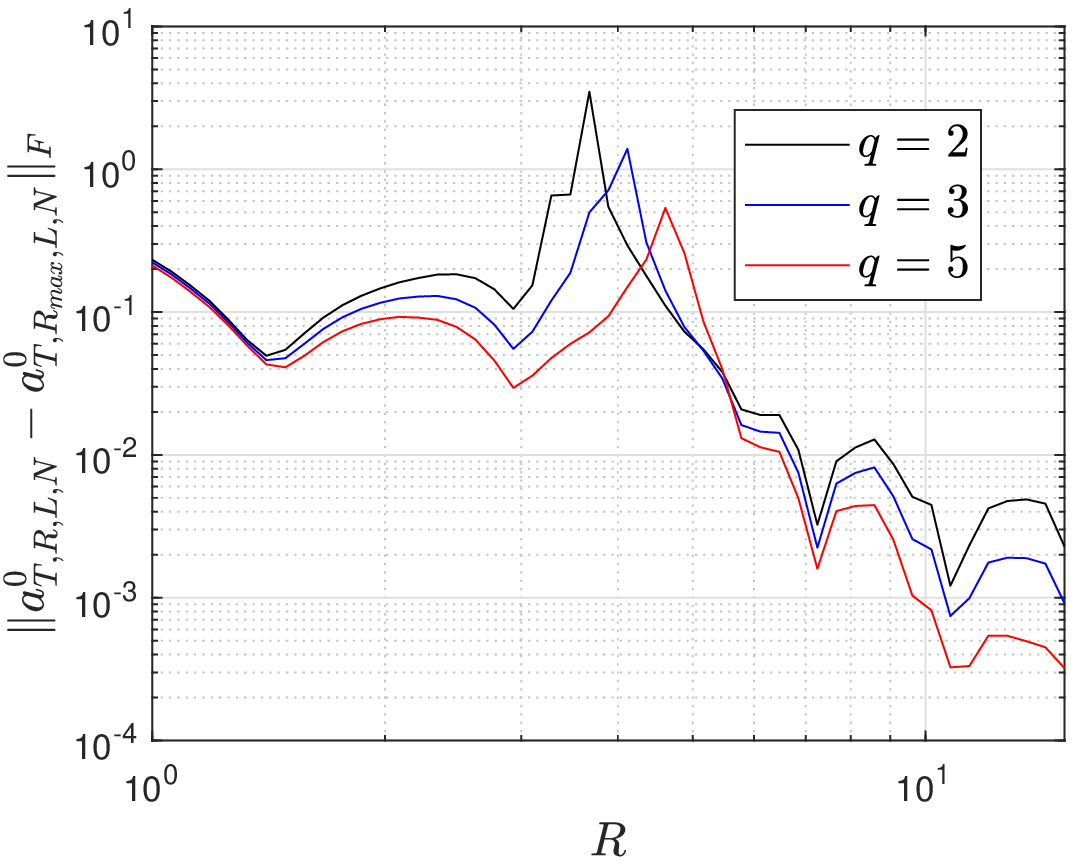}}}%
    \caption{A two dimensional random medium}%
    \label{Fig_TwoD_Random}%
\end{figure}

\textbf{Example 4. A random coefficient.} 
As yet another example of a non-periodic medium, we construct a random medium as follows: We start by choosing a large computational grid, which corresponds to a discretization of the domain $K_{R_{\max}}$ with $R_{\max} = 16$. We then  generate a sequence of uniformly distributed random variables taking values in the interval $[1,100]$, and assign these random numbers on each grid point. 
\ncorr{Next, we set a correlation length $\sigma$ (here  $\sigma = 0.25$ is chosen), and construct the random coefficient at each discretization point $\bx_i\in K_R$ (for a given $R<R_{\max}$) by taking the average of the generated random values associated to the points $\bx_j \in B_{\sigma}(\bx_i)$. We then shift the entire function to obtain a constrast ratio of $\max{a(\bx)}/\min{a(\bx)} \approx 6.78$.}
Since, the interest here is not to study the statistical error, we compute only the error
\begin{align*}
\mathcal{E}_{boundary}:= \| a^{0}_{T,R,L} - a^{0}_{T,R_{\max}, L} \|_{F},
\end{align*}
which sees the deterministic part of the overall error only; in particular the boundary error. In \cref{Fig_TwoD_Random}, the generated random coefficient along with the boundary error is depicted. All the parameter values except $h=1/40$, $\alpha \approx 19.56$, and $\beta \approx 2.88$ are the same as in example $1$. A decay for the boundary error is observed for three different choices of filters with different regularities.

\textbf{Cost comparison with the standard elliptic approach.}

\begin{figure}[ht]%
	\centering
	\subfloat[The upscaling error]{{\includegraphics[height=4.65cm]{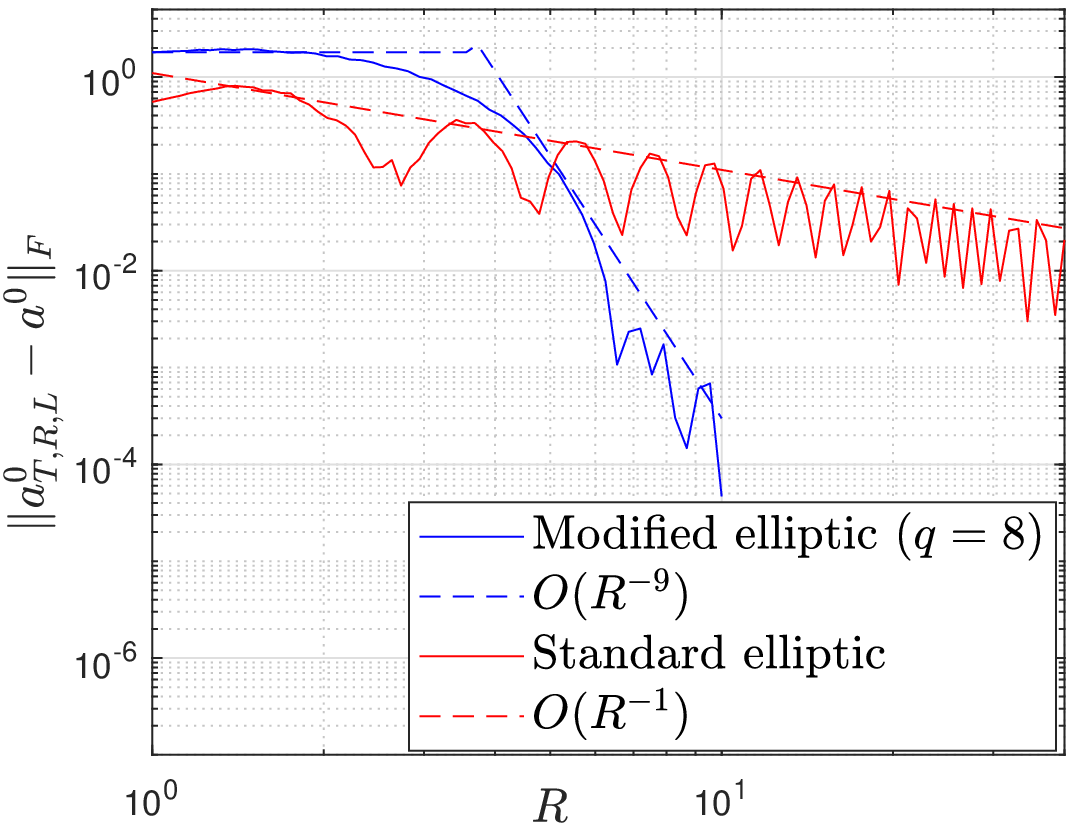} }}%
	\subfloat[Tolerance vs Computational cost]{{\includegraphics[height=4.65cm]{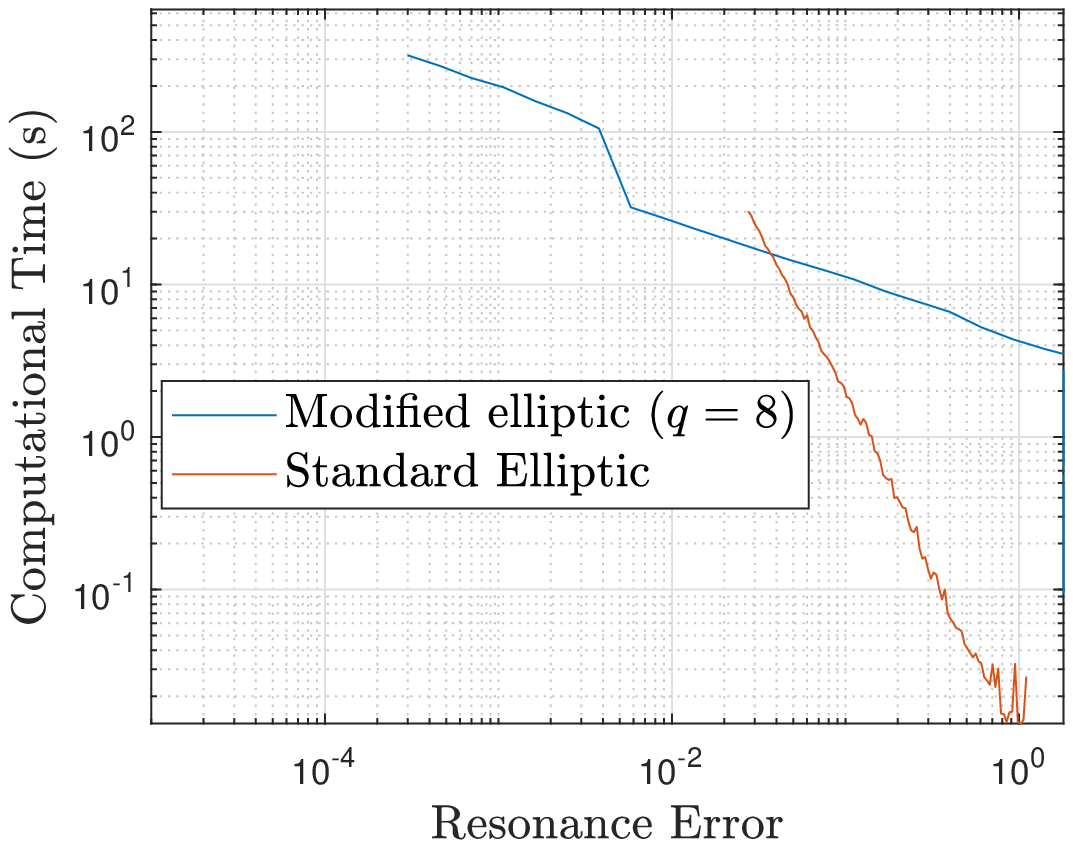}}}%
	\caption{A two dimensional periodic discontinuous medium \eqref{eqn_Coeff_Layered2D}}%
	\label{Fig_TwoD_Layered_Cost_Comp}%
\end{figure}

\ncorr{We now illustrate the elapsed computational time against tolerance, and compare the modified elliptic approach \eqref{Eqn_New_CellProblem} with the standard elliptic approach \eqref{Eqn_Approx_Periodic_Elliptic_CellProblem} approaches. The two-dimensional discontinuous coefficient from example $2$ with the contrast ratio of $10$ is used for the simulation. For both approaches the stepsize is taken to be $h=1/50$, and the computational times and the corresponding errors are recorded for a range of $R$, results of which are in Figure \ref{Fig_TwoD_Layered_Cost_Comp}. A filter with $q=8$ is used to illustrate the efficiency of the modified elliptic approach. In Figure  \ref{Fig_TwoD_Layered_Cost_Comp}(a), we compare the decay of the errors with respect to $R$. In Figure \ref{Fig_TwoD_Layered_Cost_Comp}(b), we observe that the modified elliptic approach start to become more efficient than the standard elliptic approach already for reasonably large tolerance values, $TOL  \approx 0.03$. Note that the envelopes depicted in Figure \ref{Fig_TwoD_Layered_Cost_Comp}(a) are used to obtain the results in Figure \ref{Fig_TwoD_Layered_Cost_Comp}(b). We also emphasize that, the standard elliptic approach is solved via a direct $LU$-decomposition, while the modified elliptic problem has been solved by approach 3 in subsection \ref{SubSec_CompCost}. This further explains why the modified elliptic approach is superior to the standard elliptic approach although Table \ref{Table_Cost} suggests the converse for the given parameter values. 
}

\textbf{Cost comparison with the parabolic approach from \cite{Abdulle_Arjmand_Paganoni_18b}.}

\ncorr{As discussed in the introduction, the modified elliptic approach \eqref{Eqn_New_CellProblem} has been inspired by the idea of parabolic cell-problems from \cite{Abdulle_Arjmand_Paganoni_18b}, where the following formulas were used for an approximation of the homogenized coefficient. }
\ncorr{\begin{equation*}
	\mathbf{e}_i  a^{0}_{R,L,T} \mathbf{e}_j = \int_{K_{L}} \mathbf{e}_i \cdot a(\bx) \mathbf{e}_j \mu_L(\bx)\,d\bx - 2\int_{0}^{T} \int_{ K_{L}} u^i(\bx,t) u^j(\bx,t) \mu_L(\bx) \,d\bx  \,dt,
\end{equation*} }
\ncorr{where $u^{i}$ is the solution to the parabolic equation:}
\ncorr{\begin{equation*} 
\left\lbrace
\begin{aligned}
& \frac{\partial u^i}{\partial t} - \nabla\cdot(a(\bx)\nabla u^i) = 0 & \quad &\text{in } K_{R}\times(0,T)\\
& u^i = 0 &  \quad &\text{on }  \partial K_{R} \times (0,T)\\
& u^i(\bx,0) = \nabla \cdot (a(\bx)\mathbf{e}_i) & \quad &\text{in } K_{R},
\end{aligned}
\right. 
\end{equation*}
}
\ncorr{This approximation for the homogenized coefficient results in similar exponentially decaying error bounds for the resonance error. However, the main difference is in the time dependent nature of the problem, which may pose a challenge from a computational viewpoint. An efficient numerical solver for this parabolic problem based on explicit stabilized stiff numerical solvers, such as ROCK 2, \cite{AbM01}, is presented in \cite{Abdulle_Arjmand_Paganoni_18b}. Without going into further details, we want to emphasize that ROCK 2 is an adaptive, explicit second order time-stepping method for stiff problems, with the advantage of having stability regions substantially larger than standard time integration methods. Moreover, the time-steps are computed adaptively according to a preset error tolerance value.
\begin{figure}[ht]%
	\centering
	\subfloat[The upscaling error]{{\includegraphics[height=4.65cm]{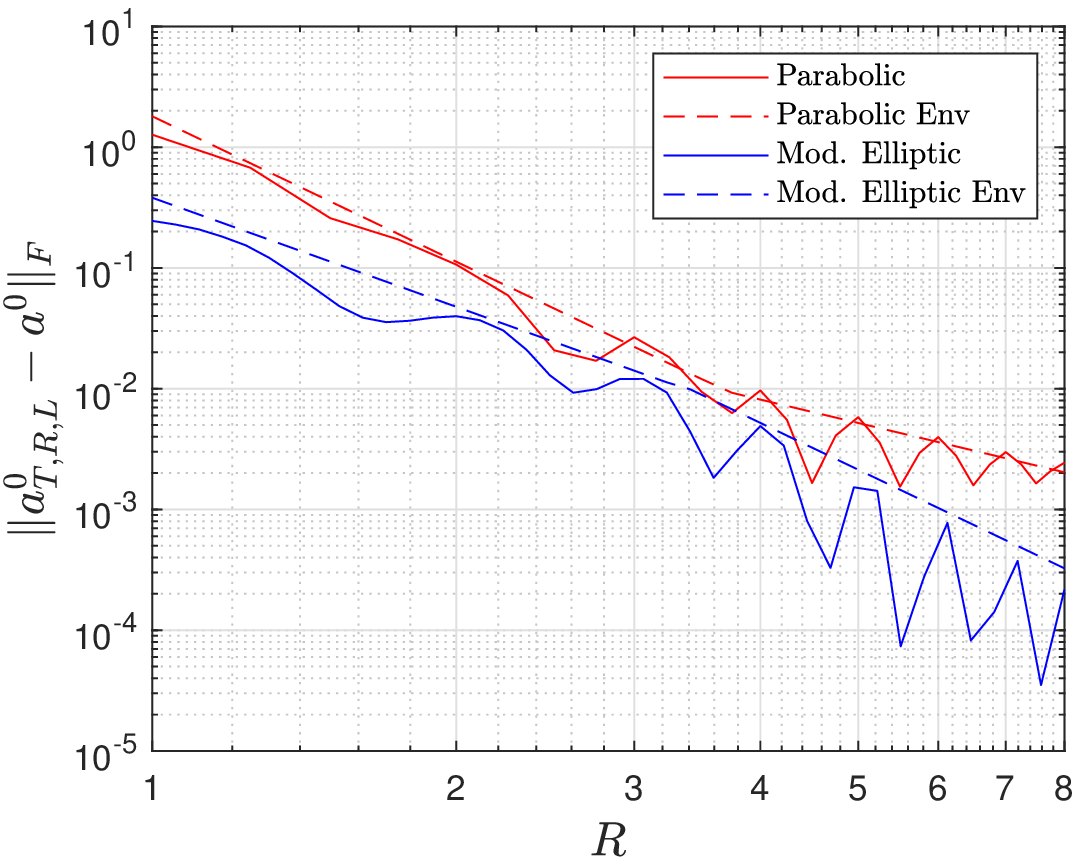} }}%
	\subfloat[Tolerance vs Computational cost]{{\includegraphics[height=4.65cm]{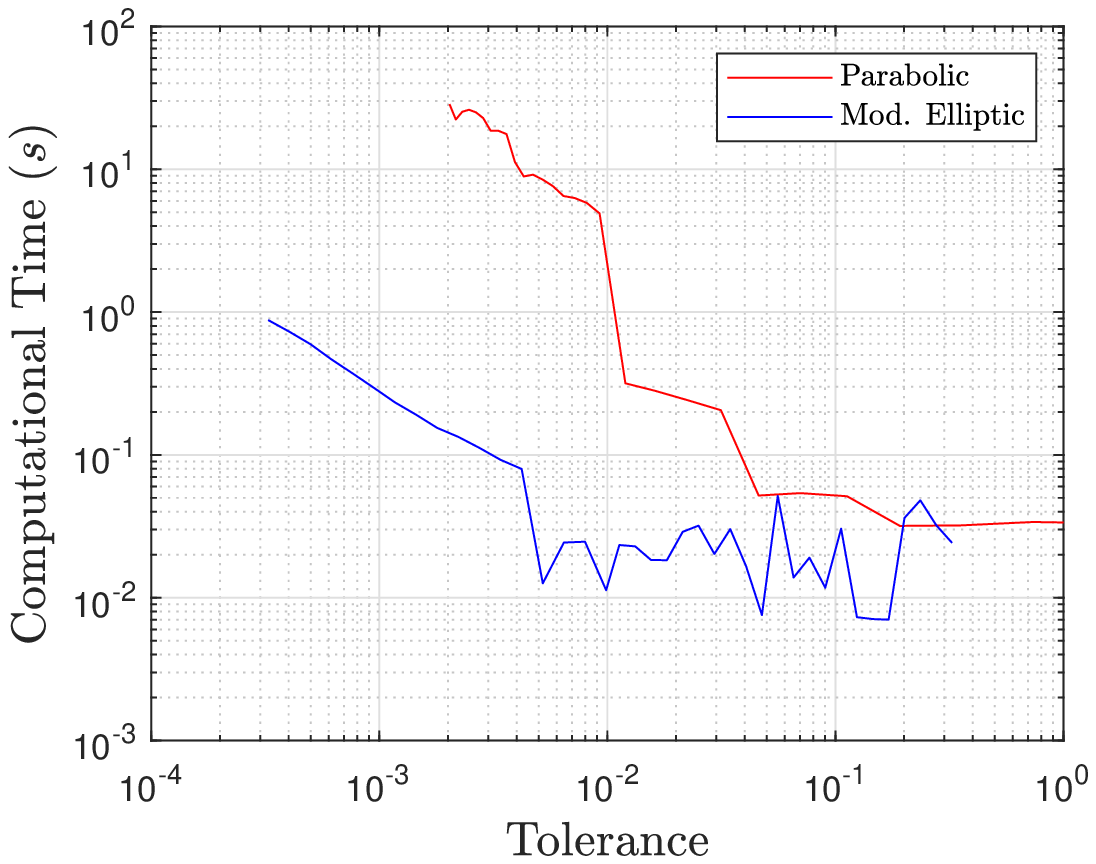}}}%
	\caption{A two dimensional periodic discontinuous medium}%
	\label{Fig_Cost_Comp_Parabolic_ModElliptic}%
\end{figure}
\begin{table}[h!]
\centering
\begin{tabular}{||c c c||} 
 \hline
  $TOL$ & $h$ & $R$ \\ [0.5ex] 
 \hline\hline
  $0.1$ & $1/8$ & $1 \leq  R \leq 1.75$ \\ 
 \hline
 $0.01$ & $1/16$ & $1.75 <  R \leq 2.5$ \\
 \hline
  $0.01$ & $1/32$ & $2.5 <  R \leq 3.5$ \\
 \hline
 $0.001$ & $1/64$ & $3.5 <  R \leq 8$ \\ [1ex] 
 \hline
\end{tabular}
\caption{Parameter values for the parabolic solver}
\label{table_Params_Parabolic}
\end{table}
\begin{table}[h!]
\centering
\begin{tabular}{||c c||} 
 \hline
  $h$ & $R$ \\ [0.5ex] 
 \hline\hline
  $1/8$ & $1 \leq  R \leq 4$ \\  
\hline
  $1/16$ & $4 <  R \leq 8$ \\ [1ex] 
 \hline
\end{tabular}
\caption{Parameter values for the modified elliptic solver}
\label{table_Params_Mod_Elliptic}
\end{table}
To compare the efficiency of the parabolic approach with the modified elliptic approach  \eqref{Eqn_New_CellProblem}, we consider again example $1$, \eqref{eqn_Coeff_Smooth2D}, with $C_1 = 2.1$ and $C_2 = 1$, and run a simulation with $q=3$ to numerically investigate the cost vs tolerance behaviour. For both simulations, the optimal values for the parameters $T$ and $L$ are chosen precisely as in example 1, and a second order method is used for the spatial discretization; the stepsize of which is denoted as $h$. For both simulations the values of $h$ are changed for increasing values of $R$ with the goal of achieving a desired tolerance value in the numerical simulation, see Tables \ref{table_Params_Parabolic} and \ref{table_Params_Mod_Elliptic}. Note that the values of $h$ are refined differently since the modified elliptic approach results in lower error tolerances already for coarse values of $h$, e.g., $h=1/8,h=1/16$, whereas in the parabolic solver, the spatial error dominates the overall error, and therefore $h$ needs to be reduced more to achieve desired error tolerances. In Figure  \ref{Fig_Cost_Comp_Parabolic_ModElliptic}, we compare (a) the actual error given by these two different approaches for increasing values of $R$, as well as (b) the elapsed computational time to reach a certain error tolerance. The envelopes for the errors depicted in Figure \ref{Fig_Cost_Comp_Parabolic_ModElliptic}(a) are used to produce the numerical results in Figure \ref{Fig_Cost_Comp_Parabolic_ModElliptic}(b). These envelopes may be seen as sharp upper bounds for the actual error, and they are used due to the non-monotonic decay of the error due to averaging. Although both methods possess exponentially decaying error bounds for the resonance error, we observe that the modified elliptic approach results in an improved cost vs error performance due to its lower discretization errors. In general, both approaches will have their advantages depending on the problem to be solved. For example, for 3-dimensional problems with possibly complicated microstructure, relying on the parabolic approach with explicit stabilized solves alleviate all linear algebra systems to be solved and the issue of preconditioning the systems.}

\section*{Acknowledgments}
The authors are grateful to Stefano Massei and Kathryn Lund for helpful discussion. This research is partially supported by Swiss National Science Foundation, grant no. 20020\_172710.

	\bibliographystyle{plain}
	\bibliography{./references}
\end{document}